\documentclass[11pt]{amsart}
\usepackage{amssymb}
\usepackage{epsfig}
\usepackage{mathdots}
 \theoremstyle{plain}
\newtheorem{theorem}{Theorem}
\newtheorem{corollary}{Corollary}

\newtheorem{proposition}{Proposition}
\theoremstyle{example}

\theoremstyle{definition}

\theoremstyle{remark}

\numberwithin{equation}{section}

\setbox0=\hbox{$+$}
\newdimen\plusheight
\plusheight=\ht0
\def\+{\;\lower\plusheight\hbox{$+$}\;}

\setbox0=\hbox{$-$}
\newdimen\minusheight
\minusheight=\ht0
\def\-{\;\lower\minusheight\hbox{$-$}\;}

\setbox0=\hbox{$\cdots$}
\newdimen\cdotsheight
\cdotsheight=\plusheight
\def\cds{\lower\cdotsheight\hbox{$\cdots$}}
\begin{document}

\title[Asymptotics of Continued Fractions]
       {Asymptotics and Sequential Closures of Continued Fractions and their Generalizations}
\dedicatory{We would like to dedicate this paper to our mathematical
father and grandfather, respectively, Basil Gordon.}
{\allowdisplaybreaks
\author{Douglas Bowman}
\address{ Northern Illinois University\\
   Mathematical Sciences\\
   DeKalb, IL 60115-2888
}
\email{bowman@math.niu.edu}

\author{James Mc Laughlin}
\address{Department of Mathematics,
West Chester University, West Chester, PA 19383}
\email{jmclaughl@wcupa.edu; http://math.wcupa.edu/\~{}mclaughlin} }
 \keywords{limit periodic continued
fractions, $q$-continued fractions, continued fractions, Poincar{\'e}-type recurrences,
$q$-series, infinite products,
asymptotics, sequential closures, Ramanujan, Cauchy distribution}
 \subjclass[2000]{Primary: 40A15, 30B70. Secondary: 39A11, 40A05, 40A20.}
\thanks{The first author's research was partially supported by NSF grant
DMS-0300126.}
\date{\today}

\begin{abstract}
Given a sequence of complex square matrices, $a_n$, consider the sequence of
their partial products, defined by $p_n=p_{n-1}a_{n}$. What can be
said about the asymptotics as $n\to\infty$ of the sequence $f(p_n)$,
where $f$ is a continuous function? A special case of our most general result 
addresses this
question under the assumption that the matrices $a_n$ are an $l_1$
perturbation of a sequence of matrices with bounded partial
products. We apply our theory to investigate the asymptotics
of the approximants of continued fractions. In particular, when a
continued fraction is $l_1$ limit  $1$-periodic of elliptic or loxodromic type, we
show that its sequence of approximants tends to a circle in $\widehat{\mathbb{C}}$,
or to a finite set
of points lying on a circle.  Our main theorem on such
continued fractions unifies the treatment of the loxodromic
and elliptic cases, which are convergent and divergent, respectively. When an
approximating sequence tends to a circle, we obtain statistical
information about the limiting distribution of the approximants. When the circle
is the real line, the points are shown to have a Cauchy distribution with
parameters given in terms of modifications of the original
continued fraction. As an example of the general theory, a detailed
study of a
$q$-continued fraction in five complex variables is provided. The most general theorem 
in the paper holds in the
context of Banach algebras. The theory is also applied to
$(r,s)$-matrix continued fractions and recurrence sequences of
Poincar\'e type and compared with closely related literature.
\end{abstract}

\maketitle

\section{Introduction}

Consider the following recurrence:

$$x_{n+1}=\frac 32 -\frac {1}{x_n}.$$
Taking $1/\infty$ to be $0$ and vice versa, then
regardless of the initial (real) value of this sequence,
it is an fact that the sequence is dense in $\mathbb{R}$. Here
is an illuminating proof.

Fix $x_0$ and note that $x_n=F^{(n)}(x_0)$, where $F$ is
the linear
fractional transformation $F(x)=(\frac 32x-1)/x$, and as usual $F^{(n)}$
denotes the $n$-th
composition of $F$ with itself. 

Next recall the correspondence
between $2 \times 2$ complex matrices and linear fractional transformations: let
a given $2\times 2$ matrix act on the complex variable $z$ by
\[
\begin{pmatrix}
a&b\\
c&d
\end{pmatrix}(z):=\frac{az+b}{cz+d},
\]
so that matrix multiplication correspond
to composition of linear fractional transformations.

A simple calculation shows that the
matrix for
$F$ is diagonalizable with diagonal elements $\alpha=3/4+i\sqrt{7}/4$ and
$\beta=3/4-i\sqrt{7}/4$. It follows that for some $A\in GL_2(\mathbb{C})$,
\[
x_n=A
\begin{pmatrix}
\alpha^n &0\\
0& \beta^n
\end{pmatrix}
A^{-1}(x_0)
=\frac{r\alpha^n+s\beta^n}{t\alpha^n+u\beta^n}=\frac{r\lambda^n+s}{t\lambda^n+u},
\]
where $\lambda=\alpha/\beta$ and $r$, $s$, $t$, and $u$ are some
complex constants. Notice that $\lambda$ is a number on
the unit circle and is not a root of unity, so that $\lambda^n$ is dense on the circle. The
conclusion follows by noting that the linear fractional transformation
\[
z\mapsto\frac{rz+s}{tz+u}
\]
is non-constant and thus a homeomorphism from
$\mathbb{\widehat{C}}$ to $\mathbb{\widehat{C}}$, and
must take the unit circle to  $\mathbb{R}\cup \{\infty\}$,
since the values of the sequence $x_n$ are real.

Letting $x_0=3/2$ and iterating the definition of $x_n$ gives that $x_n$ is the
$n$'th approximant of the continued fraction:
\begin{equation}\label{32}
3/2-\frac{1}{\displaystyle 3/2- \frac{\mathstrut 1}{\displaystyle 3/2-
\frac{\mathstrut 1}{\displaystyle 3/2-\frac{\mathstrut 1}{\ddots}}}}\,,
\end{equation}
and thus one can write down the
equality
\begin{equation}\label{32R}
\mathbb{R}=
3/2-\frac{1}{\displaystyle 3/2- \frac{\mathstrut 1}{\displaystyle 3/2-
\frac{\mathstrut 1}{\displaystyle 3/2-\frac{\mathstrut 1}{\ddots}}}}\,.
\end{equation}
This is true so long as one interprets the ``value'' of the
continued fraction as the set of limits of subsequences of its
sequence of approximants. In this paper we generalize such
equalities.

From here on we employ  space saving notations for continued
fractions.  In particular, an arbitrary finite continued fraction
\[
b_0+\cfrac{a_1}{b_1+\cfrac{a_2}{b_2+\cfrac{a_3}
{\begin{matrix}
\mathstrut\ddots& \\
 &\cfrac{a_{n-1}}{b_{n-1}+\cfrac{a_n}{b_n}}
\end{matrix}
}}}
\]
will be denoted by the expressions: $b_0+\mathop{\pmb{K}}_{i=1}^n (a_i|b_i)$,
\[
b_0+\frac{a_1}{b_1}\+\frac{a_2}{b_2}\+\cds\+\cfrac{a_n}{b_n},
\]
and
\[
b_0+\mathop{\pmb{K}}_{i=1}^n \displaystyle \frac{a_i}{b_i}.
\]
If as $n\to\infty$ this sequence of finite continued fractions converges
to a limit in $\widehat{\mathbb{C}}$ (respectively $\mathbb{C}$), then the continued fraction 
$b_0+\mathop{\pmb{K}}_{i=1}^\infty (a_i|b_i)$ is said to converge in $\widehat{\mathbb{C}}$ (respectively $\mathbb{C}$).
The limit, when it exists, is called the {\it classical limit} of the
continued fraction.

Another motivating example of our work is the following theorem, one
of the oldest in the analytic theory of continued fractions:

\begin{theorem} (Stern-Stolz, \cite{LW92,Stern,Stolz})
\label{Stern-Stolz1} Let the complex sequence $\{b_{n}\}$ satisfy $\sum
|b_{n}|$ $<\infty$. Then
\[
b_0+\mathop{\pmb{K}}_{i=1}^\infty \displaystyle \frac{1}{b_i}
\]
diverges. In fact, for $p=0,1$,
\begin{align*}
&\lim_{n \to \infty}P_{2n+p}=A_{p} \in \mathbb{C},& &\lim_{n \to \infty}Q_{2n+p}=B_{p} \in \mathbb{C},&
\end{align*}
and
\[
A_{1}B_{0}-A_{0}B_{1} = 1.
\]
\end{theorem}

Here $P_i$ and $Q_i$ are respectively the numerator and
denominator polynomials in the sequence $\{b_n\}$
obtained by simplifying the rational function one obtains by terminating
the continued fraction after $i$ terms.
The Stern-Stolz theorem shows that all continued fractions of the general form
described in the theorem tend
to two different limits, respectively $A_0/B_0$, and $A_1/B_1$. (These limits depend
on the continued fraction.) 
Although sometimes limits for continued fractions are taken to be in $\mathbb{C}$,
throughout we assume the limits for continued fractions are in $\widehat{\mathbb{C}}$.
The motivation for this is that continued fractions can be viewed as the composition
of linear fractional transformations and such functions have $\widehat{\mathbb{C}}$
as their natural domain and codomain.

Before leaving the Stern-Stolz theorem, we remark that although the theorem is
sometimes termed a ``divergence theorem'', this terminology is a bit misleading; the theorem
actually shows that although the continued fractions of this form diverge, they do so by
{\it converging}
to two limits ($A_0/B_0$ and $A_1/B_1$) according to the parity of their approximant's index.
Thus although these continued fractions diverge, they diverge in such a way that they
can still be used computationally.

A special case of the Stern-Stolz theorem is a fact about
the famous Rogers-Ramanujan continued
fraction:

\begin{equation}\label{RR}
 1+
\frac{q}{1}
\+
\frac{q^{2}}{1}
\+
\frac{q^{3}}{1}
 \+
\frac{q^{4}}{1}\+
\cds.
\end{equation}

The Stern-Stolz theorem gives that for $|q|>1$ the even and odd
approximants of \eqref{RR} tend
to two limiting functions. To see this, observe that
by the standard equivalence transformation for continued fractions, (\ref{RR})
has the same classical approximants (see the beginning of section 3 for
this terminology) as
 \[
 1+
\frac{1}{1/q}
\+
\frac{1}{1/q}
\+
\frac{1}{1/q^{2}}
 \+
\frac{1}{1/q^{2}}
\cds
\+
\frac{1}{1/q^{n}}
 \+
\frac{1}{1/q^{n}}\+
\cds
.
\]

The Stern-Stolz theorem, however, does not apply to the
following continued fraction also considered by Ramanujan:
\begin{equation}\label{R3}
\frac{-1}{1+q}
\+
\frac{-1}{1+q^2}
\+
\frac{-1}{1+q^3}
\+
\cds .
\end{equation}

Recently in  \cite{ABSYZ02} Andrews, Berndt, {\it{et al.}}  proved a claim
made
by Ramanujan in his lost notebook (\cite{S88}, p.45) about (\ref{R3}). To
describe Ramanujan's claim, we first need some notation.

Throughout take $q\in\mathbb{C}$
with $|q|<1$. The following standard notation for $q$-products
is employed:
\begin{align*}
&(a)_{0}:=(a;q)_{0}:=1,& &
(a)_{n}:=(a;q)_{n}:=\prod_{k=0}^{n-1}(1-a\,q^{k}),& & \text{ if }
n \geq 1,&
\end{align*}
and
\begin{align*}
&(a;q)_{\infty}:=\prod_{k=0}^{\infty}(1-a\,q^{k}),& & |q|<1.&
\end{align*}
Set $\omega = e^{2 \pi i/3}$. Ramanujan's claim was that, for
$|q|<1$,
{\allowdisplaybreaks
\small{
\begin{equation}\label{3lim1}
\lim_{n \to \infty}
 \left (
\frac{1}{1}
 \-
\frac{1}{1+q}
\-
\frac{1}{1+q^2}
\-
\cds
\-
\frac{1}{1+q^n +a}
\right )
=
-\omega^{2}
\left (
\frac{\Omega - \omega^{n+1}}{\Omega - \omega^{n-1}}
\right ).
\frac{(q^{2};q^{3})_{\infty}}{(q;q^{3})_{\infty}},
\end{equation}
}
}
where
{\allowdisplaybreaks
\begin{equation*}
\Omega :=\frac{1-a\omega^{2}}{1-a \omega}
\frac{(\omega^{2}q,q)_{\infty}}{(\omega q,q)_{\infty}}.
\end{equation*}
}
Ramanujan's notation is confusing, but what his claim means is
that the limit exists as $n \to \infty$ in each of the three
congruence classes modulo 3, and that the limit is given by the
expression on the right side of (\ref{3lim1}). Also, the appearance
of the variable $a$ in this formula is a bit of a red herring; from elementary
properties of continued fractions, one can derive the result for
general $a$ from information about the approximants in the $a=0$ case.

The continued fraction \eqref{32}, the Stern-Stolz theorem, and \eqref{R3} are, in fact,
examples
of the same phenomenon.
We define this phenomenon and investigate its implications.

Now \eqref{32} is different from the other two examples in that it has subsequences
of approximants tending to uncountably many limits. In a previous
paper \cite{BMcL05}, the authors presented
a theory for dealing with finite sets of limits. Here continued fractions
having finite sets of limits are unified
with examples such as \eqref{32}. Indeed,
all of the examples above, including \eqref{32}, are
special cases of a general result on continued fractions (Theorem \ref{T1} below).
More generally our results hold for a large class of sequences in Banach algebras.
To deal with
all of these cases we employ the notion of the {\it sequential closure}
of a sequence.

{\bf Definition.} The {\it sequential closure} of a sequence in a metric space
is the set of limits of convergent
subsequences.\footnote{Thus, for example, the sequence
$\{1,1,1,\dots\}$ has sequential closure $\{1\}$ although the set of limit (accumulation)
points of the set of values of the sequence is empty. Note that in a survey
paper describing some of the research in this paper, the authors previously
used the phrase ``limit set'', unaware of the use of this phrase in the theories
of discrete groups and dynamical systems. Also, the notion used here
should not be confused with the sequential closure of a set, which also
occurs in topology; thus, one should not confuse the sequential closure of a sequence
with the sequential closure of its {\it image}. We thank Peter Loeb for the suggestion of the
phrase ``sequential closure'' in the present context.
More motivation is also given in the remark following Theorem
\ref{qcfmain}.}  The sequential closure of a sequence $\{s_n\}_{n\ge
1}$ is denoted by $\copyright(s_n)$.

In this paper we study sequential closures in the specific context  of sequences of the form
\[
f\left(\prod_{i=1}^nD_i\right),
\]
where $D_i$ are elements in a unital Banach algebra and $f$ is a function
with values in a  metric space, often compact. Usually in this paper
$D_i$ is a sequence of complex matrices.

When working with infinite continued fractions we frequently omit the
symbol $\copyright$ for the sequential closure. As with \eqref{32R}, {\it equalities
involving infinite continued fractions  are
to be understood as being between the sequential closure of the continued
fraction and the set on the other side of the equation.}

This paper is divided into sections, which are now surveyed.
The main results of section 2 are Theorems \ref{DM},
\ref{DMunitconj}, and \ref{gent} which are the most general results
of the paper. Theorem \ref{DM} provides the underlying idea and is stated
at the generality of Banach algebras. In section 2 we also discuss recent
results of Beardon \cite{Bear1}, which  apply hyperbolic geometry
to study the convergence of compositions of M\"obius maps.
This approach yields short and geometrically motivated
proofs of many results in
the analytic convergence theory of continued fractions which
generalize to compositions of higher dimensional M\"obius
maps. Some of
the results of \cite{Bear1} are related to ours in as much as they
deal with generalizing the Stern-Stolz theorem. Indeed, one of the
conclusions of Theorem \ref{DM} is similar  to Theorem 4.2
from \cite{Bear1}. Overall, the main difference is that the theorems of
\cite{Bear1}, which generalize the Stern-Stolz theorem, do not
generalize the particular conclusion of the Stern-Stolz theorem that
the continued fraction's even and odd approximants tend to two
different limits, and instead focus on the fact of divergence. The
approach of this present paper is to generalize the
{\it{convergence}} of the even and odd subsequences in the Stern-Stolz theorem. 
Another difference is that the results here also give information
about the sequences $P_n$ and $Q_n$ in the Stern-Stolz 
theorem. Recovering such information from the M\"obius maps
approach seems problematic.

A special case of Theorem \ref{DM} is Theorem \ref{gent}, which is
used to prove Theorem \ref{T1}, which in turn gives detailed
information about
 the sequential closures of continued
fractions. This result is studied in detail in sections 3 and 4.
Sections 5 and 6 use Theorem \ref{gent} to study $(r,s)$-matrix continued fractions,
and linear recurrences of Poincar{\'e}  type, respectively.

Section 3 focuses on limit periodic continued fractions of elliptic
and loxodromic types. We discover a tapestry of  results which
weave together the sequential closure, modifications of the continued
fraction, and the asymptotics of the approximants of a large class
of continued fractions (including many which represent naturally
occurring special functions).   Those of elliptic type do not
converge, but we find that their sequential closures are well behaved,
and that their approximants have nice asymptotics.
Our results on the elliptic case add results to 
studies such as \cite{KL}. The foundational
result, Theorem \ref{T1}, actually treats the loxodromic and
elliptic cases simultaneously, evaluating the continued fraction
when it can have one, finitely many, or uncountably many limits
in a single formula.

Section 3 also addresses the statistics of the sequential
closure. In particular, suppose a continued fraction has an infinite
sequential closure. Then which points in the set have the ``most''
approximants tending to them, and which have the ``fewest''? Thus
for example, how can one describe the distribution of the approximants of the 
continued fraction \eqref{32R} for
$\mathbb{R}$? These questions are answered by considering the geometry of
associated linear fractional transformation. We show, for example, that
when the sequential closure is $\mathbb{R}$, the approximants follow a Cauchy distribution.
Finally, section 3
considers the converse problem and gives an algorithm for finding 
a subsequence of approximants
tending to any given point in the sequential closure.

Section 4 studies a non-trivial example of the theory.  The section concerns a particular
continued fraction with five parameters (and four degrees of freedom)
which generalizes not only the identity
\eqref{32R} above, but also the continued fraction \eqref{R3}.
In fact, the main result, Theorem \ref{qcfmain}, also gives the limit
of the Rogers-Ramanujan
continued fraction \eqref{RR} for $|q|<1$ as a special case.  Thus in one
identity, we obtain the evaluation of a continued fraction when it
has one, finitely many, or uncountably many limits.
An example of this theorem is a perturbation of \eqref{32}.
Specifically, the sequential closure of the continued fraction
\begin{equation}\label{65q}
3/2-\frac{1}{q+3/2\,-}\,\frac{1}{q^2+3/2\,-}\,\frac{1}{q^3+3/2\,-}\,\frac{1}{q^4+3/2\,-\cdots},
\end{equation}
where $|q|<1$ is complex, can be described exactly. In fact, this
sequential closure is a circle on the Riemann sphere. (Thus as a
consequence, when $|q|<1$ and $q$ is real, (\ref{65q}) always has
sequential closure $\mathbb{R}$.) Viewing this circle as a linear
fractional transformation of the unit circle
$\{z\in\mathbb{C}:|z|=1\}$,
\[
z\mapsto\frac{az+b}{cz+d},
\]
it transpires that the parameters $a$, $b$, $c$, and $d$ are special cases of the basic hypergeometric
function $_2\phi_1$, with complex base $q$.

More
generally, in \eqref{65q} if the numbers $1$ and $3/2$ are changed so that the
limiting recurrence for the convergents of the continued fraction
have distinct characteristic roots, there is a
coherent formula, in terms of the $_2\phi_1$ basic hypergeometric function,
for the sequential closure regardless of the nature of
the roots. Indeed, cases in which the
characteristic roots have a ratio that is a root of unity yield a
continued fraction with a finite set of limits. In this circumstance, the
approximants of the continued fraction converge to the set of
limits when the approximants are taken in arithmetic progressions.

{\bf Remarks and notation:}

(i) All sequential closure equalities in this paper arise  from the
situation
\[\lim_{n\to\infty}d(s_n,t_n)=0
\]
in some metric space $(X,d)$.
Accordingly,
it makes sense to define the equivalence relation $\asymp$ on sequences in $X$ by
$\{s_n\}\asymp\{t_n\}\iff \lim_{n\to\infty}d(s_n,t_n)=0$. In this situation we
refer to sequences $\{s_n\}$ and $\{t_n\}$ as being asymptotic to each other.
Abusing notation, we often write $s_n\asymp t_n$ in
place of $\{s_n\}\asymp\{t_n\}$. More generally, we frequently write sequences without braces
when it is clear from context that we are speaking of a sequence, and not the $n$th term. Identifying
the constant sequence $\{L,L,L,\dots\}$ with the value $L$, it is clear 
that the statements $\lim_{n\to\infty}x_n=L$ and $x_n\asymp L$ are equivalent.

(ii) It is a fact from general topology that given a compact topological space $X$ and
a Hausdorff space $Y$, then any continuous bijection $g:X\to Y$ must be
a homeomorphism and $g$ and its inverse must both be uniformly continuous.
Under these assumptions an immediate consequence for sequential closures is:
{\it If $\{s_n\}_{n\ge 1}$
is a sequence with values in $X$, then $\copyright(g(s_n))=g(\copyright(s_n))$}.

(iii) Another basic fact is that {\it If $\{s_n\}$ and $\{t_n\}$ are two sequences
in some metric space satisfying
 $s_n\asymp t_n$, then $\copyright(s_n)=\copyright(t_n)$.} Additionally,
if $f$ is some uniformly continuous function, then the following
sequence of implications holds: \[ s_n\asymp t_n\implies f(s_n)\asymp
f(t_n)\implies \copyright(f(s_n))=\copyright(f(t_n)).
\]

\section{Asymptotics and sequential closures of infinite products in unital Banach algebras}

 The classic theorem on the convergence of infinite products of
matrices seems to have  been given first by Wedderburn
\cite{Wedderburn1,Wedderburn2}. Wedderburn's  theorem is maybe not
as well known as it deserves to be, perhaps because Wedderburn does
not state it explicitly as a theorem, but rather gives inequalities
from which the convergence of infinite matrix products can be
deduced under an $l_1$ assumption.  Wedderburn also provides the key
inequality for establishing the invertibility of the limit, but does
not discuss this important application of his inequality. It is not
hard to see that  Wedderburn's equations  hold in any unital Banach
algebra. 
This result will be
employed in our first theorem, which gives asymptotics for
oscillatory divergent infinite products in Banach algebras. This
theorem is then applied to the Banach algebra
$\mathfrak{M}_d(\mathbb{C})$ of $d\times d$ matrices of complex
numbers topologised  using the $l_{\infty}$ norm, denoted by
$||\cdot||$.

For any unital  Banach algebra, let $I$ denote the unit. When we
use product notation for elements of a Banach algebra, or for
matrices, the product is taken from left to right; thus
{\allowdisplaybreaks
\[
\prod_{i=1}^n A_i:=A_1A_2\cdots A_n.
\]
}
Theorems with products taken in the opposite order follow from the theorems below by taking
the products in the reverse order throughout the statements and proofs.

\begin{proposition}[Wedderburn]\label{Wedderburn}
In a Banach algebra $\mathbf{U}$, let $I$ be the unit and $A_i$,  $i\ge 1$, be
a sequence.
Then $\sum_{i\ge1}||A_i||<\infty$ implies that
$\prod_{i\ge1}(I+A_i)$ converges in $\mathbf{U}$. Moreover, all the
elements of the sequence $I+A_i$ are invertible if and only if the
limit $\prod_{i\ge1}(I+A_i)$ is invertible.
\end{proposition}

The following corollary provides a convenient estimate of the
convergence rate of the product.

\begin{corollary}\label{error}
Under the conditions of Proposition \ref{Wedderburn}, let
$L=\prod_{i\ge 1}(I+A_i)$ and $P_m=\prod_{i=1}^m (I+A_i)$. Then
\begin{equation}\label{errest}
||L-P_m||\le e^{\sum_{i\ge 1}||A_i||}-e^{\sum_{1\le i \le m}||A_i||}=O\left(\sum_{i>m}||A_i||\right).
\end{equation}
\end{corollary}

\begin{proof} Put
\[
P_m=(I+A_1)(I+A_2)\cdots(I+A_m),
\]
and
\[
Q_m=(1+||A_1||)(1+||A_2||)\cdots(1+||A_m||).
\]
Expanding the product for $P_m$ gives
\begin{multline}\label{pm}
P_m=I+\sum_{1\le n_1\le m}A_{n_1}+\sum_{1 \le n_1< n_2 \le m}A_{n_1}A_{n_2}\\
+\sum_{1\le n_1< n_2<n_3\le m}A_{n_1}
A_{n_2}A_{n_3}\\
+
\cdots+\sum_{1\le n_1< n_2<\cdots<n_{m-1}\le m}A_{n_1}
A_{n_2}\cdots A_{n_{m-1}}
+A_{1}A_2\cdots A_{m}.
\end{multline}
Similarly,
\begin{multline*}
Q_m=1+\sum_{1 \le n_1 \le m}||A_{n_1}||+\sum_{1 \le n_1< n_2 \le m}||A_{n_1}||\, ||A_{n_2}||\\
+\sum_{1 \le n_1< n_2<n_3 \le m}||A_{n_1}||\,
||A_{n_2}||\, ||A_{n_3}||+
\cdots
+||A_{1}||\, ||A_2||\cdots ||A_{m}||.
\end{multline*}
Thus for $m\ge k$,
\begin{equation}\label{cauchy}
||P_m-P_k||\le (Q_m-Q_k),
\end{equation}
and 
\begin{equation}\label{1estimate}
||P_m-I||\le Q_m -1\le e^{\sum_{n\ge 1}||A_n||}-1.
\end{equation}
From the standard condition for the convergence of
infinite products of complex numbers, the convergence of $\sum_{n\ge 1}||A_n||$
implies the convergence of $\prod_{i\ge1}(1+||A_i||)$ and thus that
the sequence $Q_n$ is Cauchy. 
Hence by \eqref{cauchy}, $P_m$ is also Cauchy, and
so $\prod_{i\ge1}(I+A_i)$ exists.

Recall that an element $x$ in a Banach algebra is invertible
if $||x-I||<1$.
For $\prod_{i\ge1}(I+A_i)$ to be invertible, it is obviously necessary
that the elements of the sequence
$I+A_i$ be invertible.  We show that this is sufficient.  Since
$\sum_{i\ge1}||A_i||<\infty$, there exists $j\in\mathbb{Z}^+$
such that $\sum_{n> j}||A_n||<\log(2)$.
Then \eqref{1estimate} gives that
\[
||(I+A_{j+1})\cdots (I+A_{j+m})-I||\le e^{\sum_{n> j}||A_n||}-1.
\]
Letting $m\to\infty$ yields {\allowdisplaybreaks
\[
\lim_{m\to\infty}||(I+A_{j+1})\cdots (I+A_{j+m})-I||\le e^{\sum_{n> j}||A_n||}-1<e^{\log(2)}-1=1.
\]
}Hence $\lim_{m\to\infty}(I+A_{j+1})\cdots (I+A_{j+m})$ is
invertible. Multiplying this on the left by the invertible elements
$I+A_i$, $1\le i\le j$ gives the conclusion.
\end{proof}

\begin{proof}[Proof of Corollary]
From Proposition \ref{Wedderburn},
\begin{multline*}
||L-P_m||=\left\|\prod_{i\ge 1} (I+A_i)-\prod_{1\le i \le m} (I+A_i)\right\| \\
\le \left\|\prod_{1\le i \le m} (I+A_i)\right\|\,\left\| \prod_{ i > m} (I+A_i)-I\right\|\\
\le e^{\sum_{1 \le i \le m} ||A_i||} (e^{\sum_{ i > m}||A_i||}-1)
=O\left(\sum_{ i > m}||A_i||\right).
\end{multline*}
\end{proof}

There have  been a number of theorems  recently on the convergence
of matrix products, see \cite{Ar,Bear1,BE97,BBCM,BEN94,Hol, S00,
T97, T99}.  Also closely related to Wedderburn's theorem are
Theorems 3.7 and 3.8 of \cite{Bear1}, originally given in \cite{DT},
which gives essentially the same result, restricted to
$SL_2(\mathbf{C})$. There is also Theorem 6.1 of Borwein {\it et
al}., \cite{BBCM}, which is very similar in flavor to an earlier
result given by Artzrouni \cite{Ar}. In \cite{BL} Borwein {\it et
al.} give a generalization of the theorem in \cite{BBCM}.

Our focus here is on cases of divergence and
our results concern finding asymptotics for the $n$th partial
products. These in turn can be used to describe the sequential closures.

Here we briefly discuss the important work of Kooman,
\cite{koomR,koomM,koom}.  Kooman makes an intensive study of the
asymptotics of perturbed matrix products, recurrence sequences of
Poincar\'e type, and compositions of M\"obius transformations.
Kooman proves a number of different results on the asymptotics of
recurrences, depending on various assumptions. In particular, in
some of his results, the $l_1$ assumption is loosened, which we
maintain throughout, and other theorems treat cases of equal roots
of the characteristic equation for the recurrences. It should be
noted in the present context, however, that Kooman's results on
iterated M\"obius transformations in \cite{koom} treat the case of
the {\it outer composition sequence}, while continued fractions
arise from the {\it inner composition sequence}, and thus his work
on M\"obius transformations doesn't seem to have immediate
applications to continued fraction theory. (Recall that given a
sequence of functions $\{f_i\}$ from a set to itself, the {\it inner
composition sequence} is the sequence of functions $\{F_n\}$, where
$F_n=f_1\circ f_2\circ\cdots\circ f_n$ and the {\it outer
composition sequence}  is the sequence of functions $\{G_n\}$, where
$G_n=f_n\circ f_{n-1}\circ\cdots\circ f_1$.) Our results, which have
a stricter assumption on the perturbation in the matrix product
case, give detailed information about the relations between various
associated limits and the asymptotics for continued fractions and
generalizations. Some of this information is probably lost under
looser assumptions. But, we think that a careful study of Kooman's
results will yield new results on the convergence of continued
fractions and generalizations. Note that Kooman did use his results
on recurrences to solve a problem of Perron, see Chapter 7 of
\cite{koomR}. In section 6 we compare the applications of the
Theorem \ref{DM} below to Poincar\'e type recurrences with the
related results of Kooman and others.

We set some further conventions and fix notation. Let $G$ be a
metric space, typically a subset of ${\widehat{\mathbb{C}}}^g$,
where $\widehat{\mathbb{C}}$ is the Riemann sphere and $g$ is some
integer $g\ge1$. Here $\widehat{\mathbb{C}}$ is topologised with the
chordal metric and the corresponding product metric is employed for
${\widehat{\mathbb{C}}}^g$. (This is defined by taking the maximum
of the metrics of all the corresponding elements in two $g$-tuples.)
Let $f$ be a continuous function from a compact subset (to be
specified) of a unital Banach algebra $\mathbf{U}$, (usually
$\mathfrak{M}_d(\mathbb{C})$) to $G$. Typically we do not distinguish different
norms, the correct one being supplied from context. In a topological space
the closure of a set $S$ is denoted by $\overline{S}$.

Our first theorem is a perturbation result giving the asymptotics of
divergent infinite products in a unital Banach algebra. Although we
will only use a special case of this result, we believe the general
result is of sufficient interest to warrant inclusion, especially
since the proof of the general result requires no additional work.
We denote elements of the Banach algebra by capital letters to
suggest matrices, which is the case to which the result will be
applied.

\begin{theorem}\label{DM}
Suppose $\{M_i\}$ and $\{D_i\}$ are sequences in a unital Banach
algebra $\mathbf{U}$ such that the two sequences (for $\eta =1$
and $\eta = -1$) {\allowdisplaybreaks
\begin{equation}\label{cond1}
\left\|\left(\prod_{i=1}^nM_i\right)^{\eta}\right\|
\end{equation}
} are bounded and $\{D_i-M_i\}\in l_1({\mathbf{U}})$, that is,
\begin{equation}\label{cond2}
\sum_{i\ge 1}\|D_i-M_i\|<\infty.
\end{equation}
Let $\varepsilon_n=\sum_{i>n}||D_i-M_i||$.
Then
\begin{equation}\label{limitF}
F:=\lim_{n\to\infty}\left(\prod_{i=1}^n D_i\right)\left(\prod_{i=1}^n M_i\right)^{-1}
\end{equation}
exists and $F$ is invertible if and only if $D_i$ is invertible for all $i\geq 1$. Also,
\begin{equation}\label{err1}
\left\|F-\left(\prod_{i=1}^n D_i\right)\left(\prod_{i=1}^n M_i\right)^{-1}\right\|=O(\varepsilon_n).
\end{equation}

As sequences
\begin{equation}\label{asym}
\prod_{i=1}^n D_i\asymp F\prod_{i=1}^n M_i,
\end{equation}
and moreover
\begin{equation}\label{err2}
\left\|\prod_{i=1}^n D_i- F\prod_{i=1}^n M_i\right\|=O(\varepsilon_n).
\end{equation}
More generally, let $f$ be a  continuous function from the domain
\[
\overline{\left\{F\prod_{i=1}^n M_i:n\ge h\right\}}\cup
\bigcup_{n\ge h} \left\{\prod_{i=1}^n D_i\right\},
\]
for some integer $h\ge1$, into a metric space $G$. Then the domain
of $f$ is compact in $\mathbf{U}$ and $f(\prod_{i=1}^n D_i)\asymp
f(F\prod_{i=1}^n M_i)$. Finally
\begin{equation}
\copyright\left(\prod_{i=1}^n D_i\right)=\copyright\left(F\prod_{i=1}^n M_i\right),
\end{equation}
and
\begin{equation}
\copyright\left(f\left(\prod_{i=1}^n D_i\right)\right)=\copyright\left(f\left(F\prod_{i=1}^n M_i\right)\right).
\end{equation}
\end{theorem}

We do not assume compactness of $G$ so it is possible that
the equalities in the theorem are between empty sets. When $G$ is compact these sets
are clearly non-trivial. Note  that the conditions of the theorem imply that all the
elements $M_i$ are invertible.

An interesting special case of Theorem \ref{DM} is when the elements $M_i$ are unitary
matrices. In this paragraph the following matrix norm is used:
\[
||M||=\left(\sum_{1\le i,j\le d}|m_{i,j}|\right)^{1/2}.
\]
It is clear that $||M||=\sqrt{d}$ when $M$ is a $d\times d$ unitary matrix (for then
$||M||^2=\text{tr}(M\overline{M}^{\tt{T}})=\text{tr}(I)=d$), and thus
the hypothesis on the sequence $M_i$ is satisfied. More generally,
one can assume that the sequence of matrices $\{M_i\}$ are elements of
some subgroup of $GL_d(\mathbb{C})$ that is conjugate to the unitary group.
This case is important enough that we distinguish it in the following theorem.

\begin{theorem}\label{DMunitconj}
Let $\{M_i\}$ be a sequence of elements of a subgroup of $GL_d(\mathbb{C})$
that is conjugate to the
unitary group. Then,
if $\{D_i\}$ is a sequence $GL_d(\mathbb{C})$ and $\{D_i-M_i\}\in l_1$, all of the
conclusions of Theorem \ref{DM} hold.
\end{theorem}

The special case of Theorem \ref{DM} that will be applied in the
next section is $\mathbf{U}=\mathfrak{M}_d(\mathbb{C})$,
$M_i=M$, where $M$ be a diagonalizable complex matrix with
eigenvalues on the unit circle. Since $M$ is diagonalizable, put
$M=CEC^{-1}$. Then $M^k=CE^kC^{-1}$, and so it follows that
$\|M^k\|\le\|C\|\cdot\|C^{-1}\|$ and $\|M^k\|$ is bounded for
$k\in\mathbb{Z}$. Thus the boundedness hypothesis is satisfied and
Theorem \ref{DM} simplifies to the following.

\begin{theorem}\label{gent}
Under the above conditions,
\[
F=\lim_{n\to\infty}\left(\prod_{i=1}^n D_i\right)M^{-n}
\]
exists in $\mathfrak{M}_d(\mathbb{C})$ and $\det(F)\neq 0$ if and only if all the elements
in the sequence $\{D_i\}$ are invertible. Moreover,
$||F-(\prod_{i=1}^n D_i)M^{-n}||=O(\varepsilon_n)$,
and

(i)  $||\prod_{i=1}^n D_i- FM^n||=O(\varepsilon_n)$ . Thus $\copyright(\prod_{i=1}^n D_i)=\copyright(FM^n)$.

(ii) Let $f$ be a  continuous function from the domain
\[
\overline{\{FM^n:n\ge h\}}\cup\bigcup_{n\ge h} \left\{\prod_{i=1}^n D_i\right\},
\]
for some integer $h\ge1$, into a metric space $G$.
Then the domain of $f$ is compact and $f(\prod_{i=1}^n D_i)\asymp f(FM^n)$.
Hence  $\copyright(f\left(\prod_{i=1}^n D_i\right))=\copyright(f(FM^n))$.

\end{theorem}

Note that because $M$ is diagonalizable, this theorem enables one
to determine the exact structure of the sequential closure using Pontryagin
duality.

A special case of Theorem \ref{gent} is Proposition 1 of \cite{BMcL05} in which the eigenvalues
are roots of unity. It considers asymptotics, but not the sequential closure.
This special case is
roughly equivalent to Theorem 1.1 of \cite{short}.

\begin{proof}[Proof of Theorem \ref{DM}]
Observe that
\begin{align*}
\left(\prod_{i=1}^n D_i\right)&\left(\prod_{i=1}^n M_i\right)^{-1}
=\prod_{i=1}^n\left(\left(\prod_{j=1}^{i-1}M_j\right)D_i\prod_{j=0}^{i-1}M_{i-j}^{-1}\right)
\\
&=\prod_{i=1}^n\left(I+\left(\prod_{j=1}^{i-1}M_j\right)D_i\prod_{j=0}^{i-1}M_{i-j}^{-1}
-\left(\prod_{j=1}^{i-1}M_j\right)M_i\prod_{j=0}^{i-1}M_{i-j}^{-1}\right)
\\
&=\prod_{i=1}^n\left(I+\left(\prod_{j=1}^{i-1}M_j\right)(D_i-M_i)\prod_{j=0}^{i-1}M_{i-j}^{-1}\right)\\
&=\prod_{i=1}^n(I+A_i),
\end{align*}
where $I$ is the unit element and
\[
A_i:=\left(\prod_{j=1}^{i-1}M_j\right)(D_i-M_i)\prod_{j=0}^{i-1}M_{i-j}^{-1}.
\]
Hence
\[\|A_i\|
\leq\left\|\prod_{j=1}^{i-1}M_j\right\|\cdot\|D_i-M_i\|\cdot
\left\|\prod_{j=0}^{i-1}M_{i-j}^{-1}\right\|\leq C \|D_i-M_i\|,
\]
for some real absolute bound $C$. The second inequality followed
from the boundedness assumption on the sequences \eqref{cond1}. By
\eqref{cond2} it follows that $\sum_{i\ge0}\|A_i\|<\infty$, and so
by Proposition \ref{Wedderburn}, it follows that $F$ exists
and is invertible when the $D_i$ are invertible for $i\ge 1$.   Thus we have
proved that
\begin{equation}\label{fdif}
\lim_{n\to\infty}\left\|F-\left(\prod_{i=1}^n D_i\right)\left(\prod_{i=1}^n M_i\right)^{-1}\right\|=0.
\end{equation}
Again from the boundedness of the sequences in \eqref{cond1},
\begin{equation}\label{fdif2}
\lim_{n\to\infty}\left\|F\prod_{i=1}^n M_i-\prod_{i=1}^n D_i\right\|=0.
\end{equation}
That is,
\[
\prod_{i=1}^n D_i\asymp F\prod_{i=1}^n M_i.
\]
Using this and the boundedness of the sequences in \eqref{cond1} gives
that the domain of $f$ is compact. Thus $f$ is not only continuous, but is
uniformly continuous. This uniform continuity and \eqref{fdif2} give
\[
\lim_{n\to\infty}\left\|f\left(F\prod_{i=1}^n M_i\right)-f\left(\prod_{i=1}^n D_i\right)\right\|=0,
\]
and so
\[
f\left(F\prod_{i=1}^n M_i\right)\asymp f\left(\prod_{i=1}^n D_i\right).
\]

The sequential closure equalities in the theorem follow from the third remark in the introduction
and the error estimates follow from Corollary \ref{error} and the boundedness assumption.
\end{proof}

We conclude this section by comparing these results to some of those
from the recent paper \cite{Bear1}, which mainly focuses on applying
the hyperbolic geometry of
M\"obius maps to the convergence theory of continued fractions with
complex elements.
Consider the following two results from \cite{Bear1} that
are closely related the results of this section:

\begin{theorem}[Theorem 4.2 of \cite{Bear1}]\label{topgen}
Suppose that $G$ is a topological group whose topology is derived from
a right-invariant metric $\sigma_0$, and that $(G,\sigma_0)$ is complete. Let
$f_1,f_2,\dots$ be any sequence of elements of $G$. Then, for each $k$, there
is a neighborhood $\mathcal{N}_k$ of $f_k$ such that if, for all $j$, $g_j\in\mathcal{N}_j$,
then $(g_1\cdots g_n)(f_1\cdots f_n)^{-1}$ converges to some element $h$ of $G$.
\end{theorem}

The above theorem shares some of the structure of Theorem \ref{DM}.
In particular it gives the existence of a limit similar to the limit
$F$ in Theorem \ref{DM}.  The hypotheses are quite different,
however, and  asymptotics are not given in Theorem \ref{topgen}.
Also, sizes of the neighborhoods are not provided.

For the following corollary, some definitions involving hyperbolic
geometry are useful. A {\it M\"obius map} acting on
$\widehat{\mathbb{R}}^N$ is a finite composition of maps each of
which is an inversion or reflection in some $N-1$-dimensional
hyperplane or hypersphere in $\widehat{\mathbb{R}}^{N}$. The {\it
M\"obius group} acting on $\widehat{\mathbb{R}}^{N}$ is the group
generated by these inversions or reflections. The {\it conformal
M\"obius group}, denoted $\mathcal{M}_N$, is the subgroup of those
maps that are orientation preserving, which means that they can be
expressed as the composition of an even number of such inversions.
See \cite{Bear1, Bear2}.

\begin{corollary}[Corollary 4.3 of \cite{Bear1}] \label{dmcor}
Let $f_1,f_2,\dots$ be any sequence of M\"obius maps. Then, for each $k$, there is
a neighborhood $\mathcal{N}_k$ of $f_k$ such that if $g_j\in\mathcal{N}_j$, $j=1,2,\dots$,
then there is some M\"obius map $h$ such that for all $z$,
$\sigma(g_1\cdots g_n(z),hf_1\cdots f_n))\to 0$ as $n\to\infty$. In particular, for each point $z$,
$\lim_n g_1\cdots g_n(z)$ exists if and only if $\lim_n f_1\cdots f_n(z)$ exists.
\end{corollary}

The differences with our theorem are that the setting in Theorem
\ref{DM} is more general and  the sizes of the neighborhoods are not
given in Corollary \ref{dmcor}. However, in the case of complex
M\"obius maps, in \cite{Bear1} it is shown that the neighborhoods
$\mathcal{N}_j$ can be taken to be the set of M\"obius maps $g$ that
satisfy
\begin{equation}\label{mobcond}
||g-f_j||<\frac{1}{2^{j+2}||f_1||^2\cdots||f_{j-1}||^2||f_j||}.
\end{equation}
Here the norms are of the matrix representations of the M\"obius maps $f_i$ and $g$.

Comparing this with Theorem \ref{DM}, it can be seen that for the
case of complex M\"obius maps, unless enough of the norms $||f_i||$
are small, one expects our condition $\{D_i-M_i\}\in l_1$ to be
weaker in general, and thus our result to be stronger. Note that
Theorem \ref{DM} also gives information about the sequential closure
as well as asymptotics with error terms.  Information about the
sequential closure is implicit, however, in Corollary 4.3 of
\cite{Bear1} above.

There is another theorem in \cite{Bear1} which is also related to
Theorem \ref{DM}. In fact, it is a generalization of the Stern-Stolz
theorem presented in the introduction. Before stating the theorem, a
couple of definitions concerning  the hyperboloid model of hyperbolic
space are required.

For $x$ and $y$ in $\mathbb{R}^{N+1}$, let
\[
q(x,y)=x_1y_1+x_2y_2+\cdots+x_Ny_N-x_{N+1}y_{N+1},
\]
and
\[
\mathcal{H}_N=\{x\in\mathbb{R}^{N+1}:q(x,x)=1, \, x_{N+1}>0\}.
\]
$\mathcal{H}_{N}$ is one branch of a hyperboloid of two sheets. It
can be shown that $\mathcal{H}_N$ can be endowed with a hyperbolic
metric and that the matrix group $O^+(N+1,1)$, which preserve $q$ as
well as the condition $x_{N+1}>0$, act as isometries on this space.
Let $g$ be a M\"obius map which acts on $\mathbb{R}^N$, and hence by
the Poincar\'e extension, on $\mathbb{H}^{N+1}$. Suppose then that
$g$ corresponds to the $(N+2)\times (N+2)$ matrix $A$ which acts on
$\mathcal{H}_{N+2}$. In  \cite{Bear1}  the following beautiful
generalization of the Stern-Stolz theorem is given:

\begin{theorem}[``The General Stern-Stolz Theorem'' \cite{Bear1}] \label{gsst} Suppose
that $g_1$, $g_2$ $\dots $ are M\"obius maps in $\mathcal{M}_N$, and
that $g_n$ is represented by the $(N+2)\times (N+2)$ matrix $A_n$ as
above. If
\begin{equation}\label{SScondB}
\sum_{n=1}^{\infty}\sqrt{||A_n||^2-||I||^2}
\end{equation}
converges, then the sequence $g_1\cdots g_n$ is strongly divergent.
\end{theorem}

Consider the $N=0$ case. Then this theorem should be compared with
the case of Theorem \ref{DMunitconj} in which $H$ is unitary, and
the matrices $M_i$ represent M\"obius maps. In Theorem \ref{gsst},
\eqref{SScondB} is exactly the condition required for
$\sum_n\rho(\mathbf{j},g_n(\mathbf{j}))$ to be bounded in
$\mathbb{H}$. (Here $\rho$ is the hyperbolic metric on $\mathbb{H}$,
where $\mathbb{H}:=\{(x_1,x_2,x_3)\in\mathbb{R}^3:x_3>0\}$, and
$\mathbf{j}=(0,0,1)$.) Next, \eqref{SScondB} is sufficient to
guarantee that $\rho(\mathbf{j},g_1\cdots g_n(\mathbf{j}))$ is
finite, and thus the orbits of the product $g_1\cdots g_n$ never
leave $\mathbb{H}$. This later condition is what is meant by
``strong divergence''. Now the M\"obius maps that fix $\mathbf{j}$
are the unitary maps and $g(\mathbf{j})=\mathbf{j}$ if and only if
$||g||^2=2$. The condition \eqref{SScondB} can thus be interpreted
as saying that the elements $g_n$ approach some sequence of unitary
elements sufficiently rapidly. This is roughly the same as the
condition on the sequence $\{D_i\}$ in Theorem \ref{DMunitconj} when
$H$ is unitary. Of course the conclusion of the theorems go in
different directions.

In the next section we apply the $d=2$ case of Theorem \ref{gent} to get detailed information
about the sequential closures of continued fractions.

\section{$l_1$ limit $1$-periodic continued fractions}

The signal result of this section is Theorem \ref{T1} which for
the  unifies the evaluation (in terms of sequential closures) 
of $l_1$ limit $1$-periodic
continued fractions in loxodromic and elliptic cases. We also study
the statistics of the classical approximants in cases of sequential
closures of infinite cardinality.

First recall the correspondence between $2\times 2$ matrices and
continued fractions. It is easily understood
by remembering the
correspondence between  compositions
of linear fractional transformations and products of $2\times 2$
matrices, and then
noting that the composition of linear fractional transformations
can be written as a continued fraction. For example, observe that
for a general linear fractional transformation (avoiding cases such as $c=0$):
\begin{equation*}
\frac{az+b}{cz+d}=\frac ac+\frac{\left(\frac{bc-ad}{c^2}\right)}{\frac dc + z}.
\end{equation*}
Thus, generically,  any composition of a finite number of
non-trivial linear fractional transformations can be written as a
finite continued fraction. But to generate a continued fraction, one
does not need to work with such general linear fractional
transformations.  In particular, working with transformations of the
form
\[
\left(\frac{a_i}{b_i+z^{-1}}\right)^{-1}=\frac{b_iz+1}{a_iz}
\]
leads to
the correspondence
between matrices and continued fractions that will be used below:

{\allowdisplaybreaks
\begin{equation}\label{corr}
\left (
\begin{matrix}
& P_{n} &P_{n-1}  \\
&\phantom{as} & \phantom{as} \\
& Q_{n} & Q_{n-1}
\end{matrix}
\right )= \left (
\begin{matrix}
& b_0 & 1 \\
&\phantom{as} & \phantom{as} \\
& 1& 0\end{matrix} \right )  \left (
\begin{matrix}
& b_{1} &1  \\
&\phantom{as} & \phantom{as} \\
& a_{1} & 0
\end{matrix}
\right ) \cdots \left(
\begin{matrix}
& b_{n} &1  \\
&\phantom{as} & \phantom{as} \\
& a_{n} & 0
\end{matrix}
\right ),
\end{equation}
}
where
\begin{equation*}
\frac{P_n}{Q_n}=b_0+
\frac{a_{1}}{b_{1}} \+
 \frac{a_{2}}{b_{2}}
\+
 \frac{a_{3}}{b_{3}}
\+\,\cds \+
 \frac{a_{n}}{b_{n}}.
\end{equation*}
Here $P_n$ and $Q_n$ are the numerator and denominator
 polynomials (called {\it convergents})
in the variables $a_i$ and $b_i$ obtained
by simplifying the rational function that is the finite continued fraction. Their
ratio, $P_n/Q_n$, is called the {\it $n$th classical approximant} of the continued fraction.
When working with continued fractions with arbitrary sequences of
elements, it is assumed that at most one numerator element $a_i$ vanishes.
 Also, we assume that for
all $n$, $(a_n,b_n)\neq(0,0)$.
Moreover, if for some positive integer $N$, $a_N=0$, then
the continued fraction is finite, and thus converges. Notice that
in this case, the sequences $P_n$ and $Q_n$ still exist and are
in general
non-constant for values of $n>N$. However, in this
situation the ratio $P_n/Q_n$ will be fixed, when it is defined, and in
the cases considered, it will be proved that the ratio is defined for $n$ sufficiently large.

From \eqref{corr} one reads off immediately the fundamental recurrences for the convergents
$P_n$ and $Q_n$:
\begin{equation}\label{recur}
\left(
\begin{matrix}
P_n & P_{n-1}&\\
&\phantom{a} & \phantom{a} \\
Q_n& Q_{n-1}&
\end{matrix}
\right)
=
\left(
\begin{matrix}
P_{n-1} & P_{n-2}&\\
&\phantom{a} & \phantom{a} \\
Q_{n-1}& Q_{n-2}&
\end{matrix}
\right)
\left(
\begin{matrix}
& b_{n} &1  \\
&\phantom{as} & \phantom{as} \\
& a_{n} & 0
\end{matrix}
\right )
.
\end{equation}
Taking the determinant on both sides of \eqref{corr} gives
the {\it determinant formula} for the convergents of a
continued fraction:
\begin{equation}\label{detform}
P_nQ_{n-1}-P_{n-1}Q_n=(-1)^{n-1}a_1a_2\cdots a_n.
\end{equation}

An infinite continued fraction
\begin{equation}\label{cfgen}
\mathop{\pmb{K}}_{n=1}^{\infty}\frac{a_{n}}{b_{n}}:=\frac{a_{1}}{b_{1}} \+
 \frac{a_{2}}{b_{2}}
\+
 \frac{a_{3}}{b_{3}}
\+\,\cds
\end{equation}
is said to converge in $\mathbb{C}$ (respectively in $\mathbb{\widehat{C}}$) if
\begin{equation*}
\lim_{n \to \infty}\frac{a_{1}}{b_{1}} \+
 \frac{a_{2}}{b_{2}}
\+
 \frac{a_{3}}{b_{3}}
\+\,\cds \+
 \frac{a_{n}}{b_{n}}
\end{equation*}
exists in $\mathbb{C}$ (respectively in $\widehat{\mathbb{C}}$). Let
$\{\omega_{n}\}$ be a sequence of complex numbers. If
\begin{equation*}
\lim_{n \to \infty}\frac{a_{1}}{b_{1}} \+
 \frac{a_{2}}{b_{2}}
\+
 \frac{a_{3}}{b_{3}}
\+\,\cds \+
 \frac{a_{n}}{b_{n}+\omega_{n}}
\end{equation*}
exist, then this limit is called the {\it {modified limit of
$\mathop{\pmb{K}}_{n=1}^{\infty}(a_{n}|b_{n})$ with respect to the
sequence $\{\omega_{n}\}$}}. Detailed discussions of modified
continued fractions, as well as further pointers to the literature,
are given in \cite{LW92}. Other treatments of the analytic theory of 
continued fractions can be found in Henrici \cite{H}, Jones and Thron \cite{JT}, 
Perron \cite{P}, and Wall \cite{Wall} . Note that  by \eqref{corr} and
\eqref{recur},
\begin{equation}\label{modcf}
b_0+\frac{a_{1}}{b_{1}} \+
 \frac{a_{2}}{b_{2}}
\+
 \frac{a_{3}}{b_{3}}
\+\,\cds \+
 \frac{a_{n}}{b_{n}+\omega_{n}}=\frac{P_n+\omega_n P_{n-1}}{Q_n+\omega_n
 Q_{n-1}}.
\end{equation}

In the following theorem, the sequential closure of the sequence of
approximants of a general class of continued fractions is
computed. It transpires that the sequential closure is a circle (or a
finite subset of a circle) on the Riemann sphere. The result was obtained by
considering the problem of finding a uniform generalization of the examples
in the introduction.

Our theorem concerns the asymptotics of the continued fraction
\begin{equation}\label{basiccf}
 \frac{- \alpha\beta+q_{1}}{\alpha+\beta+p_{1}}
\+ \frac{-\alpha\beta+q_{2}}{\alpha+\beta+p_{2}} \+
 \cds
\+\frac{-\alpha\beta+q_{n}}{\alpha+\beta+p_{n}},
\end{equation}
where the sequences $p_n$ and $q_n$ approach $0$ in $l_1$ and
the constants $\alpha$ and $\beta$ are points in the complex plane. Specifically
assume that
\begin{equation}\label{pqasum}
 \sum_{n=1}^{\infty}|p_{n}|+|q_n|<\infty.
\end{equation}
Let
\begin{equation*}
\varepsilon_n:= \max\left(\sum_{i>n}|p_i|,\sum_{i>n}|q_i|\right),
\end{equation*}
and put
\[
f_n(w):= \frac{- \alpha\beta+q_{1}}{\alpha+\beta+p_{1}}
\+ \frac{-\alpha\beta+q_{2}}{\alpha+\beta+p_{2}} \+
 \cds
\+\frac{-\alpha\beta+q_{n}}{\alpha+\beta+p_{n}+w},
\]
so that $f_n:=P_n/Q_n=f_n(0)$ is the sequence of classical approximants of the
continued fraction \eqref{basiccf}.
We follow the common convention in
analysis of denoting the group of points on the unit circle by
$\mathbb{T}$, and we extend this notation in
the following way. Define the
{\it unitary characteristic} to be the map $\mathbb{T}_\lambda$  with
domain $\widehat{\mathbb{C}}$ and codomain $2^{\widehat{\mathbb{C}}}$,
given by
$\mathbb{T}_{\lambda}=\copyright(\lambda^n)\subset\widehat{\mathbb{C}}$.
From well-known facts it follows that,
\begin{equation*}
\mathbb{T}_{\lambda}=
\begin{cases}
0,  &\text{if $|\lambda|<1$;}   \\
\mathbb{T}, &\text{if $|\lambda|=1$ and $\lambda$ is not a root of unity;}      \\
\{e^{2\pi ik/m}:0\le k<m\}, &\text{if $\lambda$ is a primitive $m$th root of unity;}    \\
\infty, &\text{if $|\lambda|>1$.}
\end{cases}
\end{equation*}
We identify the unitary characteristic with its set of values.

\begin{theorem}\label{T1}
Throughout this theorem let $\{p_{n}\}_{n\geq 1}$, $\{q_{n}\}_{n \ge
1}$ be complex sequences satisfying \eqref{pqasum};
$\alpha\neq\beta$ be complex numbers with
$\alpha/\beta\in\widehat{\mathbb{C}}$. Consider the following
limits:
\begin{align}\label{abcddef}
a&=a(\alpha,\beta)=\lim_{n\to\infty}\alpha^{-n}(P_n - \beta
P_{n-1}),\\
 b&=b(\alpha,\beta)=-\lim_{n\to\infty}\beta^{-n}(P_n-\alpha P_{n-1}),\notag\\
c&=c(\alpha,\beta)=\lim_{n\to\infty}\alpha^{-n}(Q_n - \beta
Q_{n-1}),\notag\\
 d&=d(\alpha,\beta)=-\lim_{n\to\infty}\beta^{-n}(Q_n-\alpha Q_{n-1}).\notag
\end{align}
When $|\alpha|=|\beta|\neq0$, all four limits exist,
$b(\alpha,\beta)=-a(\beta,\alpha)$, and $d(\alpha,\beta)=-c(\beta,\alpha)$; the limits
for $a$ and $c$
exist if $|\alpha|>|\beta|$, while
the limits for $b$ and $d$ exist if $|\alpha|<|\beta|$.

We have
\begin{equation}\label{cflimitset}
\frac{- \alpha\beta+q_{1}}{\alpha+\beta+p_{1}}
\+ \frac{-\alpha\beta+q_{2}}{\alpha+\beta+p_{2}} \+
 \frac{-\alpha\beta+q_{3}}{\alpha+\beta+p_{3}} \+
 \cds
=\frac{a\mathbb{T}_{\alpha/\beta}+b}{c\mathbb{T}_{\alpha/\beta}+d}.
\end{equation}
(Possibly non-existent limits are always annihilated by $\mathbb{T}_{\alpha/\beta}$.)
Assuming $|\alpha|=|\beta|\neq0$, and that $q_n\neq \alpha\beta$, for $n\ge 1$,
\begin{equation}\label{cfasym}
f_n\asymp h((\alpha/\beta)^{n+1})
\,\quad\text{where}\quad\,h(z)=\frac {az+b}{cz+d};
\end{equation}
{\allowdisplaybreaks
\begin{equation}\label{detid}
\det(h)=ad-bc=(\beta-\alpha)\prod_{n=1}^\infty\left(1-\frac{q_n}{\alpha\beta}\right).
\end{equation}}
Moreover, when $|\alpha|=|\beta|\neq0$ if either $|c|\neq |d|$, or for $0\le n<m\,$,
$c\alpha^n+d\beta^n\neq 0$ when $|c|=|d|$ and $\alpha/\beta$ is a root of
unity, then as $n\to\infty$,
\begin{equation}\label{cferror}
\left|f_n-h((\alpha/\beta)^{n+1})\right|=O(\varepsilon_n).
\end{equation}
Finally, when $|\alpha|=|\beta|=1$,
the following asymptotics for the convergents $P_n$ and $Q_n$ hold as $n\to\infty$
\begin{equation}\label{Pasym}
\left|P_n-\frac{a\alpha^n+b\beta^n}{\alpha-\beta}\right|=O(\varepsilon_n)\qquad\text{and}\qquad
\left|Q_n-\frac{c\alpha^n+d\beta^n}{\alpha-\beta}\right|=O(\varepsilon_n).
\end{equation}

\end{theorem}
Note that the right hand side of \eqref{cflimitset} stands for the set
$\{(az+b)/(cz+d)|z\in\mathbb{T}_{\alpha/\beta}\}$ and reduces to a point when $\alpha/\beta$ is
$0$ or $\infty$.

{\bf Definition.} The cardinality of the sequential closure
of a continued fraction is called the \emph{rank} of the
continued fraction.
Clearly this definition makes sense for any sequence in a topological space.

Thus for the complex continued fractions covered by
Theorem \ref{T1}, the rank
belongs to the set $\mathbb{Z}^+\cup\{\mathfrak{c}\}$. (Here $\mathfrak{c}$ denotes
the cardinality of $\mathbb{R}$.)
In general, by Bernoulli's theorem on continued fractions,
complex continued fractions can also have rank
$\mathbb{\aleph}_0$.  Notice that \eqref{cflimitset}
gives a unified evaluation of $l_1$ limit $1$-periodic continued fractions in both
the loxodromic and elliptic cases.

{\bf Remark.} The identity \eqref{cflimitset} holds in the case where the continued fraction
terminates, that is when $q_n=\alpha\beta$, for some $n$,
say $n=N$.

Theorem \ref{T1} is foundational for what follows. We give two corollaries before the proofs (of Theorem \ref{T1} and its corollaries).
Further results follow the proofs. The next corollary gives enough information to
identify (up to sign)  the specific coefficients in the linear fractional transformation
$h$ in the theorem in terms of modifications of the original continued fraction. The
succeeding corollary makes that identification.

\begin{corollary}\label{modlimits}
Under the conditions of the theorem the following identities
involving modified versions of \eqref{basiccf} hold in ${\widehat{\mathbb{C}}}$.
When $|\alpha|\ge|\beta|$,
\begin{multline}\label{hinfin}
h(\infty)=\frac ac=\lim_{n\to\infty}f_n(-\beta)\\
=\lim_{n\to\infty}
\frac{- \alpha\beta+q_{1}}{\alpha+\beta+p_{1}}
\+ \frac{-\alpha\beta+q_{2}}{\alpha+\beta+p_{2}} \+
 \cds
\+\frac{-\alpha\beta+q_{n-1}}{\alpha+\beta+p_{n-1}}
\+\frac{-\alpha\beta+q_{n}}{\alpha+p_{n}}.\\
\end{multline}
When $|\alpha|\le|\beta|$,
\begin{multline}\label{h0}
h(0)=\frac bd=\lim_{n\to\infty}f_n(-\alpha)\\
=\lim_{n\to\infty}
 \frac{- \alpha\beta+q_{1}}{\alpha+\beta+p_{1}}
\+ \frac{-\alpha\beta+q_{2}}{\alpha+\beta+p_{2}} \+
 \cds
\+\frac{-\alpha\beta+q_{n-1}}{\alpha+\beta+p_{n-1}}
\+\frac{-\alpha\beta+q_{n}}{\beta+p_{n}};
\end{multline}
and for $k\in\mathbb{Z}$ when  $|\alpha|=|\beta|\neq 0$, we have
\begin{align}\notag
h((\alpha/\beta)^{k+1})&=\frac {a(\alpha/\beta)^{k+1}+b}{c(\alpha/\beta)^{k+1}+d}
=\lim_{n\to\infty}f_n(\omega_{n-k})\\\label{hvalues}&
=\lim_{n\to\infty}
\frac{- \alpha\beta+q_{1}}{\alpha+\beta+p_{1}}
\+ \frac{-\alpha\beta+q_{2}}{\alpha+\beta+p_{2}} \+
\cds
\+\frac{-\alpha\beta+q_{n}}{\alpha+\beta+p_{n}+\omega_{n-k}},
\end{align}
where
\[
\omega_{n}=-\frac{\alpha^{n} - \beta^{n}}{\alpha^{n-1} -
\beta^{n-1}}\in{\widehat{\mathbb{C}}}, \hspace{25pt} n\in\mathbb{Z}.
\]
\end{corollary}

It is also possible to derive  convergent continued fractions which
have the same limit as the modified continued fractions in Corollary
\ref{modlimits}. The key is to simply transform them via the
Bauer-Muir transformation, see \cite{LW92}.  As this will be used in
the sequel, it is presented here.

{\bf Definition} {\it The Bauer-Muir transform of a continued
fraction $b_0$ $+$\\ $\mathop{\bold{K}}(a_n|b_n)$ with respect to
the sequence $\{w_n\}$ from $\mathbb{C}$ is the continued fraction
$d_0+\mathop{\bold{K}}(c_n|d_n)$ whose canonical numerators $C_n$
and denominators $D_n$ (convergents) are given by $C_{-1}=1$,
$D_{-1}=0$, $C_n=A_n+A_{n-1}w_n$, $D_n=B_n+B_{n-1}w_n$ for
$n=0,1,2,\dots$, where $\{A_n\}$ and $\{B_n\}$ are the canonical
numerator and denominator convergents of
$b_0+\mathop{\bold{K}}(a_n|b_n)$.}

Thus the Bauer-Muir transformation gives a continued fraction whose $n$th
classical approximant is equal to the $n$th modified approximant of a given
continued fraction.

\begin{proposition}\cite{Ba,LW92,Mu}\label{BM} The Bauer-Muir transform
of $b_0+\mathop{\bold{K}}(a_n|b_n)$ with respect to $\{w_n\}$
from $\mathbb{C}$ exists if and only if
\[
\lambda_n=a_n-w_{n-1}(b_n+w_n)\neq 0
\]
for $n=1,2,3,\dots$. If it exists, then it is given by
\[
b_0+w_0+\cfrac{\lambda_1}{b_1+w_1}\+\cfrac{c_2}{d_2}\+\cfrac{c_3}{d_3}\+\cds,
\]
where $c_n=a_{n-1}s_{n-1}$, $d_n=b_n+w_n-w_{n-2}s_{n-1}$, and
$s_n=\lambda_{n+1}/\lambda_n$.
\end{proposition}

 Because of the generality of the continued fractions in
Corollary \ref{modlimits}, no substantial simplification occurs when the
Bauer-Muir transformation is applied, so we do not present
the transformed versions of the continued fractions in the corollary.

The following corollary gives (up to
a factor of $\pm 1$) the numbers $a$, $b$, $c$, and $d$
in terms
of the (convergent) modified continued fractions given in Corollary \ref{modlimits}.

\begin{corollary}\label{lftid}
When $|\alpha|=|\beta|\neq 0$
the linear fractional transformation $h(z)$ defined in Theorem
\ref{T1} has the following expression
\[
h(z)=\frac{A(C-B)z+B(A-C)}{(C-B)z+A-C},
\]
where $A=h(\infty)$, $B=h(0)$, and $C=h(1)$. Moreover, the
constants $a$, $b$, $c$, and $d$ in the theorem have the following
formulas
\[
a=sA(C-B),\quad b=sB(A-C),\quad c=s(C-B),\quad d=s(A-C),
\]
where
\[
s=\pm\sqrt{\frac{(\beta-\alpha)\prod_{n=1}^\infty\left(1-\frac{q_n}{\alpha\beta}\right)}{(A-B)(C-A)(B-C)}}.
\]
\end{corollary}

It is interesting to note that the sequence of modifications
of \eqref{basiccf} occurring in \eqref{hvalues} converge exactly
to the sequence $h((\alpha/\beta)^{n+1})$ which is asymptotic to the
approximants $f_n$ of \eqref{basiccf}.

 Dividing through the numerator and denominator of the
definition of $\omega_n$ by $\beta^{n-1}$ gives that the
 sequence $\omega_n$ occurring in \eqref{hvalues} is either a discrete
or a dense set of points on the line
\[
\frac{-\alpha\mathbb{T}+\beta}{\mathbb{T}+1},
\]
according to whether $\alpha/\beta$ is a root of unity or not. Observe
that $-\omega_{n+2}$ is the $n$th approximant
of the continued fraction
\[
\alpha+\beta+\frac {-\alpha\beta}{\alpha+\beta}\+
\frac {-\alpha\beta}{\alpha+\beta}\+
\cds\+
\frac {-\alpha\beta}{\alpha+\beta},
\]
which, except for the initial $\alpha+\beta$,
is the non-perturbed version of the continued fraction under study.
That the sequential closure of $\omega_n$ lies on a line follows from Theorem \ref{limitsetprop}
below. Combining the continued fraction for $\omega_n$ with \eqref{hvalues} and
Theorem \ref{T1} yields the intriguing
asymptotic as $k\to\infty$:
\begin{align}\notag
& \frac{- \alpha\beta+q_{1}}{\alpha+\beta+p_{1}}
\+ \frac{-\alpha\beta+q_{2}}{\alpha+\beta+p_{2}} \+
 \cds
\+\frac{-\alpha\beta+q_{k}}{\alpha+\beta+p_{k}}\\\label{newcfmod}
&\asymp
\lim_{\substack{n\to\infty}}
\frac{- \alpha\beta+q_{1}}{\alpha+\beta+p_{1}}
\+ \frac{-\alpha\beta+q_{2}}{\alpha+\beta+p_{2}} \+
\cds\notag\\\
 &\qquad\qquad\+\frac{-\alpha\beta+q_{n-1}}{\alpha+\beta+p_{n-1}} \+
\frac{-\alpha\beta+q_{n}}{p_{n}}\-
\underbrace{\frac {-\alpha\beta}{\alpha+\beta}\+
\frac {-\alpha\beta}{\alpha+\beta}\+\cds
\+
\frac {-\alpha\beta}{\alpha+\beta}}_{\text {$n-k-1$ terms}}.
\end{align}
The continued fraction on the left hand side is divergent, while its
transformed version on the right hand side asymptotically approaches
the $k$-th approximant of the continued fraction on the left as
$k\to\infty$. The relation at \eqref{newcfmod} is valid when $\sum
|p_i|+|q_i|<\infty$ and $|\alpha|=|\beta|\neq 0$, and can be viewed
as a continued fraction manifestation of Theorem \ref{gent}.

\begin{proof}[Proof of Theorem \ref{T1}]
Define
\begin{equation}\label{DMeq}
D_{n}:= \left (
\begin{matrix}
 \alpha +\beta +p_{n}& 1\\
- \alpha \beta + q_{n} &0\\
\end{matrix}
\right ), \hspace{40pt} M:= \left (
\begin{matrix}
 \alpha +\beta & 1\\
- \alpha \beta &0\\
\end{matrix}
\right ).
\end{equation}
For later use, note that
\begin{equation}\label{diageq} M=\left (
\begin{matrix}
-\beta^{-1} & -\alpha^{-1}\\
1 &1\\
\end{matrix}
\right ) \left (
\begin{matrix}
\alpha & 0\\
0 &\beta\\
\end{matrix}
\right ) \left (
\begin{matrix}
-\beta^{-1} & -\alpha^{-1}\\
1 &1\\
\end{matrix}
\right )^{-1},
\end{equation}
that for $n\in\mathbb{Z}$,
\begin{equation}\label{Mn} M^n=
\left(
\begin{array}{ll}
 \alpha ^{n+1}-\beta ^{n+1} & \alpha ^n-\beta ^n \\
 \phantom{as} & \phantom{as} \\
 -\alpha  \beta  \left(\alpha ^n-\beta ^n\right) &
   \alpha  \beta ^n-\alpha ^n \beta
\end{array}
\right)\frac{1}{\alpha - \beta },
\end{equation}
and that for $n\in\mathbb{Z}$,
\begin{equation}\label{Mnn} M^{-n}=
\left(
\begin{array}{ll}
 \alpha ^{n-1}-\beta ^{n-1} &\displaystyle{\frac{ \alpha ^{n}-\beta
 ^{n}}{\alpha \beta}}
   \\
   & \\
 \beta ^{n}-\alpha ^{n}   &\displaystyle{\frac{ \beta
 ^{n+1}-\alpha ^{n+1}}{\alpha \beta}}
\end{array}
\right)g_n,
\end{equation}
where, to save space later, we have put $g_n=(\alpha^{1-n}\beta^{1-n})/(\beta-\alpha)$.

Let $P_n$ and $Q_n$ denote the $n$th numerator and denominator
convergents of the continued fraction \eqref{basiccf}.
By the correspondence between matrices and continued fractions \eqref{corr},
{\allowdisplaybreaks
\begin{equation}\label{corr1}
\left (
\begin{matrix}
& P_{n} &P_{n-1}  \\
&\phantom{as} & \phantom{as} \\
& Q_{n} & Q_{n-1}
\end{matrix}
\right )= \left (
\begin{matrix}
& 0 & 1 \\
&\phantom{as} & \phantom{as} \\
& 1& 0\end{matrix} \right ) \prod_{j=1}^{n} D_{j}.
\end{equation}
}

Now assume $|\alpha|=|\beta|=1$.  Clearly
\[
||D_{n}-M||_{\infty} = \max \{ |p_{n}|, |q_{n}|\}.
\]
and thus
\[
\sum_{n \ge 1}||D_{n}-M||_{\infty}<\infty.
\]
It follows that the matrix $M$ and the matrices $D_{n}$ satisfy the
conditions of Theorem \ref{gent}.
Thus there exists $F\in\mathfrak{M}_2(\mathbb{C})$ defined by
 {\allowdisplaybreaks
\begin{align}\label{Feq}
&F=  \lim_{n \to \infty}\left (
\begin{matrix}
& 0 & 1 \\
&\phantom{as} & \phantom{as} \\
& 1& 0\end{matrix} \right ) \prod_{j=1}^{n} D_{j}M^{-n}\\
&= \lim_{n
\to \infty}\left (
\begin{matrix}
& P_{n} &P_{n-1}  \\
&\phantom{as} & \phantom{as} \\
& Q_{n} & Q_{n-1}
\end{matrix}
\right )M^{-n}\notag \\
& \label{expanded}
 = \lim_{n \to \infty}\left (
\begin{matrix}
& P_{n} &P_{n-1}  \\
&\phantom{as} & \phantom{as} \\
& Q_{n} & Q_{n-1}
\end{matrix}
\right )\left (
\begin{matrix}
-\beta^{-1} & -\alpha^{-1}\\
1 &1\\
\end{matrix}
\right ) \left (
\begin{matrix}
\alpha^{-n} & 0\\
0 &\beta^{-n}\\
\end{matrix}
\right ) \left (
\begin{matrix}
-\beta^{-1} & -\alpha^{-1}\\
1 &1\\
\end{matrix}
\right )^{-1} \\
&\label{prod}= \lim_{n \to \infty}\left (
\begin{matrix}
& P_{n} &P_{n-1}  \\
&\phantom{as} & \phantom{as} \\
& Q_{n} & Q_{n-1}
\end{matrix}
\right )\left(
\begin{array}{ll}
 \alpha ^{n-1}-\beta ^{n-1} &\displaystyle{\frac{ \alpha ^{n}-\beta
 ^{n}}{\alpha \beta}}
   \\
   & \\
 \beta ^{n}-\alpha ^{n}   &\displaystyle{\frac{ \beta
 ^{n+1}-\alpha ^{n+1}}{\alpha \beta}}
\end{array}
\right) \frac{\alpha^{1-n}\beta^{1-n}}{\beta-\alpha}.
\end{align}
}
Taking determinants in \eqref{expanded} gives
an expression for $det(F)$:
\[
F_{1,1}F_{2,2}-F_{1,2}F_{2,1}=\lim_{n \to
\infty}(P_{n}Q_{n-1}-P_{n-1}Q_{n})\frac{1}{(\alpha \beta)^{n}}
=-\prod_{n=1}^{\infty} \left (1-\frac{q_{n}}{\alpha \beta} \right).
\]
The last equality follows from
the determinant formula  for continued
fractions \eqref{detform}.

Consider the {\it non-terminating elliptic case}, i.e. $|\alpha|=|\beta|=1$ and $q_n\neq\alpha\beta$ for $n\ge 1$. It follows that $\det(F)\neq0$. Let
$f:GL_2(\mathbb{C})\to\widehat{\mathbb{C}}$ be given by
\[
f \left(\left(
\begin{matrix}
&u &v\\
&w &x
\end{matrix}
\right) \right)= \frac uw.
\]
Note that $f$ is continuous, and thus using Theorem \ref{gent}, is uniformly continuous
on the compact set
\[
\overline{\{FM^n:n\ge 1\}}\cup\bigcup_{n\ge 1} \left\{\left(
\begin{matrix}
&P_n &P_{n-1}\\
&Q_n &Q_{n-1}
\end{matrix}\right)\right\}.
\]
Theorem \ref{gent} and the matrix product representation of
continued fractions then give that
\[
\frac{P_n}{Q_n}\asymp f\left(FM^n\right).
\]
Hence using \eqref{Mn} and the definition of $f$,
\begin{align}\label{pnqnlim}
\frac{P_{n}}{Q_{n}} \asymp &\frac{F_{1,1}(\alpha
^{n+1}-\beta ^{n+1})+F_{1,2}(-\alpha  \beta  \left(\alpha ^n-\beta
^n\right))} {F_{2,1}(\alpha ^{n+1}-\beta ^{n+1})+F_{2,2}(-\alpha
\beta  \left(\alpha ^n-\beta ^n\right))}\\
&=\frac{(F_{1,1}- \beta
F_{1,2})\left (\frac{\alpha}{\beta}\right)^{n+1}+(\alpha F_{1,2}
-F_{1,1})} {(F_{2,1}- \beta F_{2,2})\left
(\frac{\alpha}{\beta}\right)^{n+1}+(\alpha F_{2,2} -F_{2,1})} \notag\\
&=h((\alpha/\beta)^{n+1}),\notag
\end{align}
where
\begin{equation}\label{heq}
h(z)=\frac{az+b}{cz+d},
\end{equation}
with $a=F_{1,1}- \beta F_{1,2}$, $b=\alpha F_{1,2}-F_{1,1}$,
$c=F_{2,1}- \beta F_{2,2}$,
 $d=\alpha F_{2,2}-F_{2,1}$, and $F_{i,j}\in\mathbb{C}$ are the elements of $F$.
The limit expressions for $a$, $b$, $c$, and $d$ in the theorem follow by simplifying the
constants in $h$ defined here, and then using \eqref{prod}.
Next notice that from
\eqref{prod}, the elements of the matrix $F$ are symmetric in $\alpha$ and $\beta$.
This along with the symmetry of $P_n$ and $Q_n$ as well as the definitions of
$a$, $b$, $c$, and $d$ implies that $b(\alpha,\beta)=-a(\beta,\alpha)$ and also
that $d(\alpha,\beta)=-c(\beta,\alpha)$. The limits \eqref{abcddef} are clearly
invariant of the size of $|\alpha|$ (since $P_n$ is a polynomial of degree $n$
in $\alpha$), so they all exist under just the
assumption $|\alpha|=|\beta|\neq0$.  Note that we can
compactly express the definition of $a$, $b$, $c$, and $d$ in the following matrix equation:
\[
\left(
\begin{matrix}
a& b& \\
c& d&
\end{matrix}
\right)
=\left(
\begin{matrix}
F_{1,1}& F_{1,2}& \\
F_{2,1}& F_{2,2}&
\end{matrix}
\right)
\left(
\begin{matrix}
1& -1&\\
-\beta& \alpha&
\end{matrix}
\right).
\]
The product formula
for $ad-bc$ follow immediately by taking the determinant
and using the expression for $\det(F)$ above.
Solving for $F$ gives
\begin{equation}\label{Fmatrix}
F=
\left(
\begin{matrix}
a& b&\\
c& d&
\end{matrix}
\right)
\left(
\begin{matrix}
\alpha& 1&\\
\beta& 1&
\end{matrix}
\right)
\frac 1{\alpha-\beta}.
\end{equation}

Now $h:\widehat{\mathbb{C}}\to\widehat{\mathbb{C}}$
is a continuous bijection when $\det(h)\neq0$. Put $\lambda=\alpha/\beta$.
From \eqref{pnqnlim} and the remarks in the introduction,
\[
\copyright\left(\frac{P_n}{Q_n}\right)=\copyright(h(\lambda^{n+1}))
=h(\copyright(\lambda^{n+1}))=h(\mathbb{T}_\lambda),
\]
and so \eqref{cflimitset} is proved.

To prove the asymptotics for $P_n$ and $Q_n$ employ Theorem \ref{gent} (i)
to obtain
\begin{equation}\label{PP}
\left(
\begin{matrix}
P_n& P_{n-1}&\\
Q_n& Q_{n-1}&
\end{matrix}
\right)
\asymp
FM^n.
\end{equation}
Substituting \eqref{Mn} and \eqref{Fmatrix} into \eqref{PP} yields
\begin{align*}\allowdisplaybreaks
&\left(
\begin{matrix}
P_n& P_{n-1}&\\
Q_n& Q_{n-1}&
\end{matrix}
\right)\\
&\qquad\asymp
\left(
\begin{matrix}
a& b&\\
c& d&
\end{matrix}
\right)
\left(
\begin{matrix}
\alpha& 1&\\
\beta& 1&
\end{matrix}
\right)
\left(
\begin{matrix}
\alpha^{n+1}-\beta^{n+1}& \alpha^{n}-\beta^{n}&\\
-\alpha\beta (\alpha^{n}-\beta^{n})& \alpha\beta^n-\beta\alpha^n&
\end{matrix}
\right)
\frac 1{(\alpha-\beta)^2}\\
&\qquad=
\left(
\begin{matrix}
a& b&\\
c& d&
\end{matrix}
\right)
\left(
\begin{matrix}
\alpha^{n+1}& *&\\
\beta^{n+1}& *&
\end{matrix}
\right)
\frac 1{\alpha-\beta}\\
&\qquad=
\left(
\begin{matrix}
a\alpha^{n+1}+b\beta^{n+1}& *&\\
c\alpha^{n+1}+d\beta^{n+1}& *&
\end{matrix}
\right)
\frac 1{\alpha-\beta}.
\end{align*}
Thus the sequences $P_n$ and $Q_n$ have the claimed asymptotics by Theorem \ref{gent}.

Put $A_n=a\alpha^n+b\beta^n$ and $B_n=c\alpha^n+d\beta^n$ and observe that
\begin{align*}
\left| f_n-h(\lambda^{n+1})\right|&=
\left|\frac {P_n}{Q_n}-\frac {A_n}{B_n}\right|
\le \left|\frac{P_nB_n-A_nB_n}{Q_nB_n}\right|+\left|\frac{A_nB_n-Q_nA_n}{Q_nB_n}\right|\\
&\le\left|\frac{1}{Q_n}\right|\varepsilon_n+\left|\frac{A_n}{Q_nB_n}\right|\varepsilon_n,
\end{align*}
and this error is $O(\varepsilon_n)$ providing that $B_n$ is bounded away
from $0$. (Recall that
$Q_n\asymp B_n/(\alpha-\beta)$.) It is easy to see that $B_n$ is bounded away from $0$ under
exactly the two conditions given in the theorem.
The restriction $|\alpha|=|\beta|=1$ can be loosened to $|\alpha|=|\beta|\neq0$
by employing the equivalence transformation in which the
numerator elements of the continued fraction are divided through by $|\alpha|^2$, and the
denominator elements by $|\alpha|$.

Now consider the {\it non-terminating loxodromic case }, i.e.  $|\alpha|\neq |\beta|$ and $q_n\neq \alpha\beta$, and 
assume that $\alpha\beta\neq 0$.
Since the continued fraction \eqref{basiccf} is
of loxodromic type, it converges by Theorem 28, p.p. 151--152 of \cite{LW92}. By Theorem 5.1
of \cite{koomR} the recurrence for the convergents,
\[
Y_n=(\alpha+\beta+p_n)Y_{n-1}+(-\alpha\beta+q_n)Y_{n-2}
\]
has a basis of solutions $\{u_n(\alpha,\beta),v_n(\alpha,\beta)\}$ satisfying
\begin{equation}\label{convasym}
\lim_{n\to\infty}\frac{u_n(\alpha,\beta)}{\alpha^n}
=\lim_{n\to\infty}\frac{v_n(\alpha,\beta)}{\beta^n}=1.
\end{equation}
Thus there exist $r(\alpha,\beta)$ and $s(\alpha,\beta)$ such
that
\[
P_n=r(\alpha,\beta)u_n(\alpha,\beta)+s(\alpha,\beta)v_n(\alpha,\beta).
\]
By \eqref{convasym} when $|\alpha|>|\beta|$, $
\lim_{n\to\infty}\alpha^{-n}{P_n}=r(\alpha,\beta) $; similarly
$\lim_{n\to\infty}\alpha^{-n}$ $Q_n$ exists. Hence the limits for
the constants $a$ and $c$ in  \eqref{abcddef} exist. Since $|\alpha|>|\beta|$,
$\mathbb{T}_{\alpha/\beta}=\infty$. Thus the right
hand side of \eqref{cflimitset} reduces to $a/c$. Now,
\[
\frac ac=\lim_{n\to\infty}\frac{P_n-\beta P_{n-1}}{Q_n-\beta Q_{n-1}}=\lim_{n\to\infty}f_n(-\beta).
\]
It is well-known, see \cite{LW92} p. 160, that $\lim_{n\to\infty} f_n(-\beta)=\mu$
when $|\alpha|>|\beta|>0$, where $\mu$ is the classical limit of the continued fraction. 
The case $|\beta|>|\alpha|>0$ follows by symmetry.

Next consider the {\it non-terminating loxodromic case} where $|\alpha|>\beta=0$, with $q_n\neq 0$ for $n\ge 1$.
By taking an equivalence
transformation, this is equivalent to the case $\alpha=1$, $\beta=0$. Again by Theorem 5.1
of \cite{koomR} the recurrence for the convergents,
\[
Y_n=(1+p_n)Y_{n-1}+q_nY_{n-2}
\]
has a solution $\{u_n(\alpha,\beta)\}$ satisfying $\lim_{n\to\infty}u_n=1$. By the Poincar\'e--Perron
theorem (see the beginning of section 6) there also exists a solution $\{v_n(\alpha,\beta)\}$ satisfying $\lim_{n\to\infty} v_{n+1}/v_n=0$. This implies that there exists a $k>0$ such that
$v_n\neq 0$ for $n\ge k$. Select such a $k$. Then for $n>k$,
\[
\lim_{n\to\infty}\frac{v_{n}}{v_k}=
\lim_{n\to\infty} \frac{v_{k+1}}{v_k}\cdot\frac{v_{k+2}}{v_{k+1}}\cdot\cdots\cdot\frac{v_{n}}{v_{n-1}}
=0,
\]
and thus it can be concluded that $\lim_{n\to\infty} v_n=0$.
Since $P_n$ and $Q_n$ are linear combinations of $u_n$ and $v_n$,
it follows that the limits for $a$ and $c$ in \eqref{abcddef} reduce to
$a=\lim_{n\to\infty} P_n$ and $c=\lim_{n\to\infty}Q_n$.
It is known that when $\alpha=1$ and $\beta=0$, the continued fraction
converges in $\widehat{\mathbb{C}}$ to $\lim_{n\to\infty} P_n/Q_n=a/c$;
see \cite{LW92}, p.p. 151--152. When $\alpha=1$ and $\beta=0$ the right hand side of \eqref{cflimitset}
simplifies thus:
\[
\frac{a\mathbb{T}_{\alpha/\beta}+b}{c\mathbb{T}_{\alpha/\beta}+d}
=\frac{a{\infty}+b}{c{\infty}+d}=\frac ac,
\]
which as has just been shown is the classical limit of
the continued fraction.

Now consider the {\it terminating elliptic case}, that is assume
that $|\alpha|=|\beta|=1$, and
that $q_N=\alpha\beta$ for some $N\ge 1$. Theorem \ref{gent} still gives that
$F$ exists and that the asymptotics \eqref{Pasym} also hold.
It follows that $a$, $b$, $c$, and $d$ exist, although $ad-bc=0$. It is easy
to check that in this situation $a=a^*\alpha^{-N}P_{N-1}$, $b=b^*\beta^{-N}P_{N-1}$,
$c=a^*\alpha^{-N}Q_{N-1}$, and $d=b^*\beta^{-N}Q_{N-1}$, where the numbers
$a^*$ and $b^*$ are defined by the limits
\[
a^*=\lim_{k\to\infty}\alpha^{-k}(P^*_k-\beta P^*_{k-1}),
\]
\begin{equation}\label{b*}
b^*=-\lim_{k\to\infty}\beta^{-k}(P^*_k-\alpha P^*_{k-1}),
\end{equation}
and $P^*_k$ is the $k$th numerator convergent of the continued fraction
\begin{equation}\label{cftail}
b_N+\mathop{\pmb{K}}_{i=1}^{\infty}\frac{-\alpha\beta+q_{N+i}}{\alpha+\beta+p_{N+i}}.
\end{equation}
(Note that we have use the easily proved
identities $P_{N+k}=P^*_kP_{N-1}$ and $Q_{N+k}=P^*_kQ_{N-1}$ which follow
from the assumption that $q_N=\alpha\beta$.)
Because \eqref{cftail} has no vanishing numerators, it follows that it
is impossible for both $a^*$ and $b^*$ to vanish. Indeed, what we have already
proved above applies to the continued fraction \eqref{cftail} 
(so by \eqref{detid} $a^*d^*-b^*c^*\neq 0$).
Computing
the right hand side of \eqref{cflimitset} gives:
\begin{equation*}
\begin{split}
\frac{a\mathbb{T}_{\alpha/\beta}+b}{c\mathbb{T}_{\alpha/\beta}+d}&=
\frac{a^*\alpha^{-N}P_{N-1}\mathbb{T}_{\alpha/\beta}+b^*\beta^{-N}P_{N-1}}
{a^*\alpha^{-N}Q_{N-1}\mathbb{T}_{\alpha/\beta}+b^*\beta^{-N}Q_{N-1}}\\
&=
\left(\frac{P_{N-1}}{Q_{N-1}}\right)
\frac{a^*\alpha^{-N}\mathbb{T}_{\alpha/\beta}+b^*\beta^{-N}}
{a^*\alpha^{-N}\mathbb{T}_{\alpha/\beta}+b^*\beta^{-N}}
=\frac{P_{N-1}}{Q_{N-1}},
\end{split}
\end{equation*}
since at least one of the limits $a^*$ or $b^*$ is non-zero.
The condition that $|\alpha|=|\beta|=1$ can now be
loosened to just $|\alpha|=|\beta|\neq0$ by applying an equivalence
transformation.

Finally, consider the {\it terminating loxodromic case}, i.e. $|\alpha|>|\beta|$ and
$q_N=\alpha\beta$. As $\mathbb{T}_{\alpha/\beta}=\infty\in\widehat{\mathbb{C}}$,
we need to show that the limits $a$ and $c$ exist and that $a/c=P_{N-1}/Q_{N-1}$.
The proofs in the non-terminating loxodromic cases apply here as well 
and show that the limits $a$ and $c$ exist.
Moreover, $a=a^*\alpha^{-N}P_{N-1}$, and
$c=a^*\alpha^{-N}Q_{N-1}$, where $a^*$
is as defined in the terminating elliptic case above.
Again, the same argument as in the non-terminating loxodromic cases give
that $\lim_{k\to\infty}\alpha^{-k}P^*_k$ exists.
Moreover, from the Poincar\'e-Perron theorem on recurrences, see section 6,
$\lim_{k\to\infty}P_{k+1}^*/P^*_k$ exists,
and is equal to either $\alpha$ or $\beta$. Thus, since
$\alpha,\beta\in\mathbb{C}$, there follows $P^*_k\neq 0$ for $k$
sufficiently large. Hence it follows that $P_{N-1}/Q_{N-1}=a/c$ in this
case. The case $|\alpha|<|\beta|$ and $q_N=\alpha\beta$ follows by symmetry.
\end{proof}

\begin{proof}[Proof of Corollary \ref{modlimits}]
 \eqref{hinfin} and \eqref{h0} follow immediately from
the value of a modified continued fraction \eqref{modcf}, with
$\omega_{n} = -\beta$ and $\omega_{n} =-\alpha$, respectively,
and the limit expressions for $a$, $b$, $c$, and $d$.

Let $f$ denote the function from the last proof. To get \eqref{hvalues}, observe that
{\allowdisplaybreaks
\begin{align*}
&h(\lambda^{k+1})=f(FM^k)=f\left(\lim_{n\to\infty} \left(
\begin{matrix}
&P_n &P_{n-1}\\
&Q_n &Q_{n-1}
\end{matrix}\right)M^{-n}M^k\right)\\
&=f\left(\lim_{n\to\infty}
\left(
\begin{matrix}
&P_n &P_{n-1}\\
&Q_n &Q_{n-1}
\end{matrix}\right)M^{-(n-k)}\right)\\
&=f\left(\lim_{n\to\infty} \left(
\begin{matrix}
P_n &P_{n-1}\\
Q_n &Q_{n-1}
\end{matrix}\right)
\left(
\begin{matrix}
 \alpha ^{n-k-1}-\beta ^{n-k-1} &\displaystyle{\frac{ \alpha ^{n-k}-\beta
 ^{n-k}}{\alpha \beta}}
   \\
   & \\
 \beta ^{n-k}-\alpha ^{n-k}   &\displaystyle{\frac{ \beta
 ^{n-k+1}-\alpha ^{n-k+1}}{\alpha \beta}}
\end{matrix}
\right)g_{n-k} \right)\\
&=\lim_{n\to\infty}\frac{(\alpha^{n-k-1}-\beta^{n-k-1})P_n-(\alpha^{n-k}-\beta^{n-k})P_{n-1}}
{(\alpha^{n-k-1}-\beta^{n-k-1})Q_n-(\alpha^{n-k}-\beta^{n-k})Q_{n-1}}\\
&=\lim_{n\to\infty}\frac{P_n-
\displaystyle{\frac{\alpha^{n-k}-\beta^{n-k}}{\alpha^{n-k-1}-\beta^{n-k-1}}}P_{n-1}}
{Q_n-\displaystyle{\frac{\alpha^{n-k}-\beta^{n-k}}{\alpha^{n-k-1}-\beta^{n-k-1}}}Q_{n-1}}\\
&=\lim_{n\to\infty}\frac{P_n+\omega_{n-k}P_{n-1}}
{Q_n+\omega_{n-k}Q_{n-1}},
\end{align*}
} where
\[
\omega_j:=-\frac{\alpha^{j}-\beta^{j}}{\alpha^{j-1}-\beta^{j-1}}.
\]
The result now follows from \eqref{modcf}.
\end{proof}

\begin{proof}[Proof of Corollary \ref{lftid}]
The expression for $h(z)$ follows immediately using algebra from
\eqref{hinfin}, \eqref{h0}, and \eqref{hvalues} with $k=-1$. The expressions for $a$,
$b$, $c$, and $d$ follow by using \eqref{detid} along with the fact that the coefficients
in the two expressions for the
linear fractional transformation must be equal up to a constant factor.
\end{proof}

Note that putting $k=0$ and $k=-1$ in \eqref{pnqnlim}  gives the following identities:
\begin{align}\label{k0-1}
h(\lambda)&=\frac{F_{1,1}}{F_{2,1}},\\
h(1)&=\frac{F_{1,2}}{F_{2,2}}. \notag
\end{align}

One naturally wonders just how effectively the parameters $a$, $b$, $c$, and $d$ in
Theorem \ref{T1}
can be computed. In the next section,  a particular continued fraction is considered
which generalizes one of Ramanujan's, as well as \eqref{32}, and  these
parameters explicitly are computed as
well-behaved meromorphic functions of the variables in the continued fraction.
Thus, for the $q$-continued fraction studied in the next section,
the parameters can not only be computed, but also have
nice formulas.

An interesting special case of Theorem \ref{T1} occurs when
$\alpha$ and $\beta$ are distinct $m$-th roots of unity ($m \geq
2$). In this situation the continued fraction
\[
\frac{- \alpha\beta+q_{1}}{\alpha+\beta+p_{1}} \+
\frac{-\alpha\beta+q_{2}}{\alpha+\beta+p_{2}} \+
\frac{-\alpha\beta+q_{3}}{\alpha+\beta+p_{3}} \+
\frac{-\alpha\beta+q_4}{\alpha+\beta+p_{4}} \+ \cds\]
becomes limit
periodic and the sequences of approximants  in the $m$ different
arithmetic progressions modulo $m$ converge. The corollary below,
which is also proved in \cite{BMcL05}, is an easy consequence of
Theorem \ref{T1}. Note that by Theorem \ref{T1} \eqref{basiccf} can also have a finite
sequential closure in the more general case that $\alpha/\beta$
is root of unity, a case not covered in the following corollary.

\begin{corollary}\label{c1}
Let $\{p_{n}\}_{n\geq 1}$, $\{q_{n}\}_{n \ge 1}$ be complex
sequences satisfying
\begin{align*}
 &\sum_{n=1}^{\infty}|p_{n}|<\infty,& &\sum_{n=1}^{\infty}|q_{n}|<\infty.&
\end{align*}
Let $\alpha$  and $\beta$ be distinct roots of unity and let $m$ be
the least positive integer such that $\alpha^m=\beta^{m}=1$ . Define
{\allowdisplaybreaks
\begin{equation*}
G:= \frac{- \alpha\beta+q_{1}}{\alpha+\beta+p_{1}} \+
\frac{-\alpha\beta+q_{2}}{\alpha+\beta+p_{2}} \+
\frac{-\alpha\beta+q_{3}}{\alpha+\beta+p_{3}} \+ \cds.
\end{equation*}
} Let $\{P_{n}/Q_{n}\}_{n=1}^{\infty}$ denote the sequence of
approximants of $G$. If $q_{n}\not = \alpha\beta$ for any $n \geq
1$, then $G$ does not converge. However, the sequences of numerators
and denominators in each of the $m$ arithmetic progressions modulo
$m$ do converge. More precisely,
 there exist complex numbers $A_{0}, \dots  , A_{m-1}$ and $B_{0}, \dots  , B_{m-1}$
such that, for $0 \leq i< m$, {\allowdisplaybreaks
\begin{align}\label{ABlim}
&\lim_{k \to \infty} P_{m\,k+i}=A_{i}, & &\lim_{k \to \infty}
Q_{m\,k+i}=B_{i}.&
\end{align}
Extend the sequences $\{A_i\}$ and $\{B_i\}$ over all integers by
making them periodic modulo $m$ so that (\ref{ABlim}) continues to
hold. Then for integers $i$,
\begin{equation}\label{lol}
A_i=\left(\frac{A_1-\beta A_0}{\alpha-\beta}\right)\alpha^i
+\left(\frac{\alpha A_0-A_1}{\alpha-\beta}\right)\beta^i,
\end{equation}
and
\begin{equation}\label{lolb}
B_i=\left(\frac{B_1-\beta B_0}{\alpha-\beta}\right)\alpha^i
+\left(\frac{\alpha B_0-B_1}{\alpha-\beta}\right)\beta^i.
\end{equation}

Moreover,
\begin{equation}\label{detAB}
A_iB_{j}-A_{j}B_i=-(\alpha\beta)^{j+1}
\frac{\alpha^{i-j}-\beta^{i-j}}{\alpha-\beta}
\prod_{n=1}^{\infty}\left(1-\frac{q_{n}}{\alpha\beta}\right).
\end{equation}
} Put  $\alpha := \exp( 2 \pi i a/m)$,  $\beta := \exp (2 \pi i
b/m)$, $0 \leq a <b <m$, and $r:=m/\gcd(b-a,m)$. Then $G$ has rank $r$
and its sequential closure is the finite set
in $\widehat{\mathbb{C}}$ given by $\{A_j/B_j:1\le j\le r\}$.
Finally, for $k \geq 0$ and $1 \leq j \leq r$,
{\allowdisplaybreaks
\begin{equation*}
\frac{A_{j+kr}}{B_{j+kr}}=\frac{A_{j}}{B_{j}}.
\end{equation*}
}
\end{corollary}

\begin{proof}
Let $M$ be as in Theorem \ref{T1}. It follows from \eqref{diageq}
that
\begin{equation}\label{mpower}
M^{j}= \left (
\begin{matrix}
\displaystyle{
 \frac{{{\alpha}}^{1 + j} - {{\beta}}^{1 + j}}
   {{\alpha} - {\beta}} }
& \displaystyle{ \frac{{{\alpha}}^j - {{\beta}}^j}
   {{\alpha} - {\beta}}} \\
\phantom{as} & \phantom{as} \\
- \displaystyle{ \frac{{\alpha}\,{\beta}\, \left( {{\alpha}}^j -
{{\beta}}^j \right) }{{\alpha} - {\beta}}}

& \displaystyle{ \frac{- {{\alpha}}^j\,{\beta}   +
     {\alpha}\,{{\beta}}^j}{{\alpha} - {\beta}} }
\end{matrix}
\right ),
\end{equation}
and thus that {\allowdisplaybreaks
\begin{align*}
&M^{m}=\left (
\begin{matrix}
&1 &0 \\
&0 &1
\end{matrix}
\right ),& &M^{j} \not = \left (
\begin{matrix}
&1 &0 \\
&0 &1
\end{matrix}
\right ),& &1 \leq j < m.&
\end{align*}
}

Let the matrix $F$ be as in Theorem \ref{T1}. From
 \eqref{PP},
 {\allowdisplaybreaks
\begin{equation}\label{abproof}
\lim_{n \to \infty} \left (
\begin{matrix}
& P_{mn+i} &P_{mn+i-1} \\
&\phantom{as} & \phantom{as} \\
& Q_{mn+i} & Q_{mn+i-1}
\end{matrix}
\right ) = \lim_{n \to \infty} F\,M^{mn+i}=F\,M^{i}.
\end{equation}
} This proves \eqref{ABlim}.

Now let $A_i:=\lim_{n\to\infty}P_{mn+i}$, and
$B_i:=\lim_{n\to\infty}Q_{mn+i}$. Notice by definition that the
sequences $\{A_i\}$ and $\{B_i\}$ are periodic modulo $m$. It easily
follows from  (\ref{abproof}) that {\allowdisplaybreaks
\begin{equation*}
\left (
\begin{matrix}
& A_{i} &A_{i-1} \\
&\phantom{as} & \phantom{as} \\
& B_{i} & B_{i-1}
\end{matrix}
\right ) = \left (
\begin{matrix}
& A_{j} &A_{j-1} \\
&\phantom{as} & \phantom{as} \\
& B_{j} & B_{j-1}
\end{matrix}
\right ) M^{i-j}.
\end{equation*}
} (\ref{mpower}) also gives that {\allowdisplaybreaks
\begin{equation}\label{2da}
A_{i} = A_{j} \displaystyle{
 \frac{{{\alpha}}^{1 + i-j} - {{\beta}}^{1 +i- j}}
   {{\alpha} - {\beta}} }
-A_{j-1} \displaystyle{ \frac{{\alpha}\,{\beta}\,
       \left( {{\alpha}}^{i-j} - {{\beta}}^{i-j} \right) }{{\alpha} - {\beta}}},
\end{equation}
and
\begin{equation}\label{2db}
B_{i} = B_{j} \displaystyle{
 \frac{{{\alpha}}^{1 + i-j} - {{\beta}}^{1 +i- j}}
   {{\alpha} - {\beta}} }
-B_{j-1} \displaystyle{ \frac{{\alpha}\,{\beta}\,
       \left( {{\alpha}}^{i-j} - {{\beta}}^{i-j} \right) }{{\alpha} - {\beta}}}.
\end{equation}
} Thus
\begin{equation*}
A_{i}B_{j}-A_{j}B_{i}=
 \frac{\left( {A_j}\,{B_{-1 + j}} -
        {A_{-1 + j}}\,{B_j} \right) \,{\alpha}\,{\beta}\,
      \left( {{\alpha}}^{i - j} - {{\beta}}^{i - j} \right) }
      {{\alpha} - {\beta}}.
\end{equation*}
Equations (\ref{lol}) and (\ref{lolb}) follow from (\ref{2da}) and
(\ref{2db}) by setting $j=1$. (\ref{detAB}) follows after applying
the determinant formula
\begin{align*}
A_{j}B_{j-1}-A_{j-1}B_{j}
&= -\lim_{k \to \infty} \prod_{n=1}^{mk+j}(\alpha\beta-q_{n})\\
&= - (\alpha\beta)^{j}\prod_{n=1}^{\infty} \left (
1-\frac{q_{n}}{\alpha\beta} \right ).
\end{align*}
Since $\sum_{j=1}^{\infty}|q_{j}|$ converges to a finite value,
 the infinite product on the right side converges.

For the continued fraction to converge,
$A_{i}B_{i-1}-A_{i-1}B_{i}=0$ is required. However, (\ref{detAB})
shows that this is not the case.
\end{proof}

\subsection{Distribution of approximants}

Let $\mathbb{T}^\prime$ denote the image of $\mathbb{T}$ under $h$,
that is, the sequential closure of the sequence $\{f_n\}$. The asymptotic for
$f_n$ given in Theorem \ref{T1} is
{\allowdisplaybreaks\begin{equation}\label{fn} f_n\asymp
h(\lambda^{n+1}),
\end{equation}
} where $h$ is the linear fractional transformation defined in the
theorem and $\lambda=\alpha/\beta$.

Some observations can immediately be  made. It is well
known that when $\lambda$ is not a root of unity, $\lambda^{n+1}$ is uniformly
distributed on
$\mathbb{T}$.  However, the linear fractional transformation $h$
stretches and compresses  arcs of the circle $\mathbb{T}$, so that the distribution of
 $h(\lambda^{n+1})$ in arcs of $\mathbb{T}^\prime$ is no longer uniform.
(Recall uniform distribution on a curve happens when as
$n\to\infty$  each segment of the curve get's the proportion of
the first $n$ points equal to the ratio of the segment's length to
the length of the whole curve.)   Additionally, $\mathbb{T}^\prime$ may
not be compact in $\mathbb{C}$. So we consider a probability measure on
$\mathbb{T}^\prime$ giving the probability of an element
$h(\lambda^{n+1})$ being contained in a subset of $\mathbb{T}^\prime$.
This measure is easy to write down.  Let
$S\subset\mathbb{T}^\prime$, then $h^{-1}(S)$ is a subset of the
unit circle. Then since $\lambda^n$ is uniformly distributed on
$\mathbb{T}$, $P(S):=\mu(h^{-1}(S))/2\pi$ gives the probability
that for any $n$, $h(\lambda^n)\in S$. Here $\mu$ denotes the
Lebesgue measure on $\mathbb{T}$.
Note that $P$ depends entirely on $h$, and thus only on the
parameters $a$, $b$, $c$, and $d$.

In general $f_n\notin \mathbb{T}^\prime$, but because of
\eqref{fn}, as $n\to\infty$, the terms of the sequence $f_n$ get closer and
closer to the sequence $h(\lambda^{n+1})$ which lies on $\mathbb{T}^\prime$.  Thus
we speak of $P$ as the {\it {limiting probability measure for the sequence $f_n$ with
respect to $\mathbb{T}^\prime$.}} When the sequential closure is $\mathbb{R}$, we speak
of the {\it limiting probability density function}.

More specifically, \eqref{fn} implies that there is a one-to-one correspondence
between the convergent subsequences of $h(\lambda^{n+1})$ and those of $f_n$ such
that the corresponding subsequences tend to the same limit. As $h$ is a
homeomorphism and $\lambda^n$ is uniformly distributed on $\mathbb{T}$,
 it follows that the probability of an element of $\copyright(f_n)$ being
contained in a subset $S$ of $\mathbb{T}^\prime$ is exactly $P(S)=\mu(h^{-1}(S))/2\pi$.

Fortunately, this distribution is completely controlled by the
known parameters $a$, $b$, $c$, and $d$. The following theorem gives
the points on the sequential closures whose neighborhood arcs have the greatest
and least
concentrations
of approximants.

\begin{theorem}\label{limitsetprop}
When $\alpha/\beta=\lambda\in\mathbb{T}$ is not a
root of unity and $cd\neq 0$, the points on
\[\frac{a\mathbb{T}+b}{c\mathbb{T}+d}
\]
with the highest and lowest concentrations of approximants are
{\allowdisplaybreaks
\[
\frac{\displaystyle{\frac{a}{c}|c|+\frac{b}{d}|d|}}{|c|+|d|}
\qquad\text{and}\qquad\frac{-\displaystyle{\frac{a}{c}|c|+\frac{b}{d}|d|}}{-|c|+|d|},
\]
} respectively. If either $c=0$ or $d=0$, then all points on the
sequential closure have the same concentration. The radius of the sequential closure
circle in $\mathbb{C}$ is {\allowdisplaybreaks\[
\left|\frac{\alpha-\beta}{|c|^2-|d|^2}
\prod_{n=1}^{\infty}\left(1-\frac{q_n}{\alpha\beta}\right)\right|,
\]
} and its center is the complex point
\[
\frac{|h(1)|^2(h(-1)-h(i))+|h(-1)|^2(h(i)-h(1))+|h(i)|^2(h(1)-h(-1))}
{h(1)(\overline{h(i)}-\overline{h(-1)})+
h(-1)(\overline{h(1)}-\overline{h(i)})+h(i)(\overline{h(-1)}-\overline{h(1)})}.
\]

The sequential closure is a line in $\mathbb{C}$ if and only if $|c|=|d|$, and in
this case the point of least concentration is $\infty$.
\end{theorem}

\begin{proof}Let $g(\theta)=h(e^{i\theta})$. Thus $g(\theta)$ parametrizes
$\mathbb{T}^\prime$ for $\theta\in [0,2\pi]$ and $e^{i\theta}$
moves with a uniform speed around $\mathbb{T}$ as $\theta$ moves
uniformly from $0$ to $2\pi$. Then $g(\theta)$ moves around
$\mathbb{T}^\prime$ at different speeds depending on how the
length $g(\theta)$ change with $\theta$. Accordingly, we wish to
compute the rate of change of the length of $g(\theta)$ with
respect to $\theta$. We then wish to know when this value is
minimum and maximum. To this end put
\[
l(\theta):=\int_0^\theta|g^{\prime}(t)|dt.
\]
Accordingly, $l^{\prime}(\theta)=|g^{\prime}(\theta)|$. An easy computation gives
\[
l^{\prime}(\theta)=\frac{|ad-bc|}{|c|^2+|d|^2+c\overline{d}e^{i\theta}+\overline{c}de^{-i\theta}},
\]
and thus
\[
l^{\prime\prime}(\theta)=i\frac{|ad-bc|(\overline{c}de^{-i\theta}-c\overline{d}e^{i\theta})}
{(|c|^2+|d|^2+c\overline{d}e^{i\theta}+\overline{c}de^{-i\theta})^2}.
\]
Clearly $l^{\prime\prime}(\theta)=0$ if and only if
$e^{i\theta}=\pm |c|d/c|d|$. Plugging these values into $h$ gives
the points where the length of $g(\theta)$ is changing most and
least with respect to $\theta$.

To find the radius of $\mathbb{T}^\prime$, one computes $l(2\pi)/2\pi$:
\begin{align*}
\frac{l(2\pi)}{2\pi}
&=\frac{1}{2\pi}\int_0^{2\pi}\frac{|ad-bc|}{|c|^2+|d|^2+c\overline{d}e^{i\theta}+\overline{c}de^{-i\theta}}d\theta\\
&=\frac{|ad-bc|}{2\pi i}\oint\frac{dz}{(c+dz)(\overline{d}+\overline{c}z)},
\end{align*}
where the contour on the last integral is the unit circle. A
routine evaluation by the residue theorem along with \eqref{detid}
gives the result.
The center can
 easily be computed as it is the circumcenter of the triangle formed by
any three points on the circle, for example, $z_1=h(1)$, $z_2=h(-1)$,
and $z_3=h(i)$. The  well-known formula for the circumcenter of three non-collinear points
in the complex plane
\[
\frac{|z_1|^2(z_2-z_3)+|z_2|^2(z_3-z_1)+|z_3|^2(z_1-z_2)}
{z_1(\overline{z_3}-\overline{z_2})+
z_2(\overline{z_1}-\overline{z_3})+z_3(\overline{z_2}-\overline{z_1})}
\]
thus gives the center of the sequential closure circle.
The final conclusions of the theorem follow immediately from the formulas for
the points of highest and lowest concentration.
\end{proof}

\begin{corollary}\label{c3}
If the sequential closure of the continued fraction
in \eqref{basiccf} is a line in $\mathbb{C}$, then
the point of highest concentration of approximants in the sequential closure
is exactly
\[
x_0=\frac{h(\infty)+h(0)}{2}=\frac 12\left(\frac ac +\frac bd \right),
\]
the average of the first two modifications of \eqref{basiccf} given in Corollary \ref{modlimits}.
Moreover, if the sequential closure is $\mathbb{R}$, then the limiting probability density
function for the approximants is given by the Cauchy density function
\begin{equation}\label{pdf}
p(x)=\frac{\delta}{\pi[(x-x_0)^2+\delta^2]},
\end{equation}
where $\delta$ is the dispersion  (scale) parameter given by
\[
\delta=\frac{h(\infty)-h(0)}{2i}=\frac{1}{2i}\left(\frac ac -\frac bd \right).
\]
\end{corollary}

For real argument period $1$ continued fractions (not limit periodic) the density function has
been studied before (somewhat more informally than here); see \cite{OOT1,OOT2}.

\begin{proof}
If the sequential closure is a line, then Theorem \ref{limitsetprop} implies that $|c|=|d|$.
The same theorem also implies that the point of highest concentration
is given by
\[\frac{\frac{a}{c}|c|+\frac{b}{d}|d|}{|c|+|d|}.
\]
When $|c|=|d|$, this simplifies to
\[
\frac12 \left(\frac ac +\frac bd\right),
\]
the average of $h(\infty)$ and $h(0)$.

Suppose the sequential closure is $\mathbb{R}$. Let the point $x\in\mathbb{R}$
be related to the point $z$ on the unit circle via
\[
x=h(z)=\frac{az+b}{cz+d},
\]
and suppose $z=e^{i\theta}$.  Let $\theta_0\in (0,2\pi]$ be
the angle for which $z$ is mapped to $\infty$ by $h(z)$, and put $z_0=e^{i\theta_0}$.
Let $p(x)$ denote the probability
density function and  let $f_i$ denote the $i$-th
approximant of \eqref{basiccf}. Then for any
interval $[a,b]$,
\begin{align*}
\int_{a}^{b} p(x) d x &= \lim_{n \to \infty} \frac{\# \{f_i  \in[a,b]\}_{ 0 \leq i \leq n}}{n}\\
& =\frac{\mu (h^{-1}([a,b]))}{2 \pi}, \notag
\end{align*}
where,  the
second equality follows from remarks made in the discussion
preceding Theorem \ref{limitsetprop}. In particular,
\[
\int_{-\infty}^{x} p(t) d t = \frac{\text{length of the arc
clockwise from $z_0$ to $z$}}{2 \pi} = \frac{\theta_0 - \theta}{2
\pi}.
\]
Using the Fundamental Theorem of Calculus, one obtains
\begin{equation*}
\begin{split}
p(x) &= \frac{-1}{2 \pi} \frac{d \theta}{d x} = \frac{-1}{2 \pi i z}
\frac{d z}{d x}=\frac{ad-bc}{2 \pi i(cx-a)(dx-b)}\\
&=\frac{h(\infty)-h(0)}{2\pi i(x-h(\infty))(x-h(0))}=\frac{\delta}{\pi[(x-x_0)^2+\delta^2]},
\end{split}
\end{equation*}
where the parameters are as given in the corollary.
\end{proof}

Note that our arguments above (in the discussion preceeding Theorem \ref{limitsetprop}, 
and the proofs of Theorem \ref{limitsetprop} and its corollary) apply
to any uniformly distributed sequence on $\mathbb{T}$ and thus the 
first part of the following corollary follows.

\begin{corollary} Let $\{s_n\}$ be a sequence that is uniformly 
distributed on $\mathbb{T}$ and let $h(z)=(az+b)/(cz+d)$ be a 
linear fractional transformation that maps $\mathbb{T}$ to $\mathbb{R}$. Then
the sequence $\{h(s_n)\}$ has a Cauchy distribution on $\mathbb{R}$ with 
parameters $x_0=(a/c+b/d)/2$ and $\delta=(a/c-b/d)/(2i)$. 

Conversely, every Cauchy distribution on $\mathbb{R}$ arises as such a M\"obius transformation
of a uniformly distributed sequence on $\mathbb{T}$.

\end{corollary}

\begin{proof} Because of the above discussions we need only to prove
the converse direction. Let a sequence $\{s_n\}$ have a Cauchy distribution 
on $\mathbb{R}$ with parameters
$x_0$ and $\delta$. Then the sequence $\{s_n/\delta -x_0/\delta\}$ has
a standard Cauchy distribution centered at $0$ with dispersion parameter $1$ on $\mathbb{R}$. 
Then the transformation $(iz+1)/(-iz+1)$ maps this sequence to a uniformly 
distributed sequence on $\mathbb{T}$, see \cite{Mc}. Finally note that the 
composition of these two maps of the sequence $\{s_n\}$ is a M\"obius transformation.
\end{proof}

The circular Cauchy distributions
of \cite{Mc} are special cases of the family of distribution functions on circles in $\hat{\mathbb{C}}$
obtained by applying an arbitrary non-trivial M\"obius transformation to a uniform
distribution on the unit circle. We call this family the M\"obius-Cauchy distribution family.
It is clearly closed under the full M\"obius
group, unlike the circular Cauchy distributions which are closed under only those
M\"obius maps that fix the unit circle, \cite{Mc}. As (real) Cauchy distributions are closed under
the real M\"obius group, the M\"obius-Cauchy family forms a natural generalization
of the usual real Cauchy distribution that also includes the family of circular Cauchy distributions
as a special case.

\subsection{Computing subsequences of approximants converging to any point on
the sequential closure}

We recall one of the main conclusions of Theorem \ref{T1}. Namely,
that if $\sum |p_n|<\infty$, $\sum |q_n|<\infty$,
$|\alpha|=|\beta|=1$ and $\lambda=\alpha/\beta$ is not a root of
unity, then $f_n$, the $n$-th approximant of $\mathop{\pmb{K}}(-\alpha  \beta +q_i | \alpha
+\beta+p_i)$,  satisfies
\[
f_n \asymp h(\lambda^{n+1}):=\frac{a \lambda^{n+1} +b}{c \lambda^{n+1}
+d},
\]
for some $a$, $b$, $c$ and $d \in \mathbb{C}$. Thus the approximants
densely approach a circle in the complex plane and a natural
question is the following: is it possible to explicitly determine a
subsequence of approximants converging to $h(e^{2 \pi i \theta})$,
for any $\theta \in [0,1)$? The following algorithm solves this problem 
assuming knowledge of the regular continued fraction for $\theta$.

Let $\lambda =e^{2 \pi i \gamma}$, $\gamma \in (0,1)$ and let
$\{a_n/b_n\}$ denote the sequence of even indexed
approximants in the regular
continued fraction expansion of $\gamma$. Since $\lambda$ is
not a root of unity, it follows that $\gamma$ is irrational.  For
real $z$, let $\langle z\rangle$ denote the fractional part of $z$. Thus
$\langle z\rangle=z-\lfloor z\rfloor$.
Let $\theta \in [0,1)$
and, for $n \geq 1$, let $r_n$ denote the least positive integer
satisfying $0\le r_n/b_n-\theta<1/b_n$. For any positive integer
$x$,
\[
x \gamma - \theta = x\left(\gamma -\frac{a_n}{b_n}\right)
+\frac{x a_n -r_n}{b_n} + \left (\frac{r_n}{b_n} - \theta \right).
\]

Since $\gcd(a_n,b_n)=1$, there exists a non-negative integer $x<b_n$
satisfying $a_nx\equiv r_n \pmod {b_n}$. Let $k_n$ be this solution.
Since $(a_nk_n-r_n)/b_n\in\mathbb{Z}$, it follows that
\[
\langle k_n\lambda-\theta\rangle=\left\langle k_n\left(\gamma -\frac{a_n}{b_n}\right)
+ \left (\frac{r_n}{b_n} - \theta \right)\right\rangle.
\]

If the sequence $\{ k_n\}$ is unbounded, let $\{j_n\}$ be a strictly
increasing subsequence. If $\{ k_n\}$ is bounded, replace each $k_n$
by $k_n+b_n$ and once again let $\{j_n\}$ be a strictly increasing
subsequence. From the theory of regular continued fractions we have
that, in either case, \[k_n\left|\gamma
-\frac{a_n}{b_n}\right|<(k_n+b_n)\left|\gamma
-\frac{a_n}{b_n}\right|<\frac{2}{b_n},\] and thus that
\[
\langle j_n \gamma - \theta\rangle \to 0.
\]
 It now follows
that $f_{j_n -1} \asymp h(\gamma^{j_n}) \to h(e^{2 \pi i \theta})$.
Thus
\[
\lim_{n \to \infty} f_{j_n -1} =h(e^{2 \pi i \theta}).
\]

Note that for rational $\lambda=m/n$, one takes approximants
in arithmetic progressions modulo $n$ to obtain the subsequences
tending to the discrete sequential closure.

\subsection{Related work}

We are aware of four other places where work related to the results of this section
was given previously.
Two of these were motivated by the identity \eqref{3lim1} of Ramanujan. The
first paper is \cite{ABSYZ02} which gave the first proof of \eqref{3lim1}.
The proof in \cite{ABSYZ02} is particular to the continued fraction \eqref{3lim1}.
However, section 3 of \cite{ABSYZ02} studied the recurrence $x_n=(1+a_{n-1})x_{n-1}
-x_{n-2}$ and showed that when $\sum_{n\ge1}|a_n|<\infty$, the sequence
$\{x_n\}$ has six limit points and that moreover a continued fraction whose
convergents satisfies this recurrence under the $l_1$ assumption tends
to three limit points (Theorem 3.3 of  \cite{ABSYZ02}). The paper does not
consider other numbers of limits, however.  Moreover, the role of the sixth roots of
unity in the recurrence is not revealed. In the section 6
of the present paper, we treat the general case in which recurrences can have
a finite or uncountable number of limits. Previously in \cite{BMcL05} we treated such
recurrences with a finite number $m>1$ of limits as well as the associated
continued fractions.

Ismail and Stanton \cite{IS05} gave a proof of \eqref{3lim1} and
also obtained Corollary \ref{ramgentheor} below. Their method was to
use properties of orthogonal polynomials and they obtained theorems
on continued fractions with any finite number of limits.
Unfortunately the method has limitations on the perturbing
sequences, requiring them to be real. On the other hand, when the
approach applies, it yields a formula for the limits of the
continued fraction in terms of associated orthogonality measures.
The general theorem of \cite{IS05}, Theorem 5.2, is actually a
simple application of  Theorem 40, of Nevai \cite{Nev}. The overall
approach of \cite{IS05} was actually employed previously in
\cite{AI}. \cite{IS05} also contains a number of other beautiful
explicit new continued fraction evaluations, similar to
\eqref{3lim1}.

We also mention here, in as much as it deals with the convergence 
of subsequences of approximants
of continued fractions, the results obtained when the approximants 
of a continued fractions form normal families. In these cases there 
are theorems for expressing the limiting function (of a convergent 
subsequence of approximants) as 
Stieltjes integrals. See for example Henrici \cite{H} and Wall \cite{Wall}. 
The work of the present paper deals, however, with the pointwise 
limits of the continued fractions, rather than the limits of 
subsequences of functions of the variable $z$ in a continued fraction.

Finally, we briefly compare our results with a theorem
of Scott and Wall \cite{SW,Wall}.
Consider the continued fraction
\begin{equation}\label{swt}
\frac 1{b_1}\+\frac 1{b_2}\+\frac 1{b_3}\+\cds.
\end{equation}
\begin{theorem}[Scott and Wall]\label{swa}
If the series $\sum|b_{2p+1}|$ and $\sum|b_{2p+1}s_p^2|$, where $s_p=b_2+b_4+\cdots+b_{2p}$,
converge, and $\liminf |s_p|<\infty$, then the continued fraction \eqref{swt} diverges.
The sequence of its odd numerator and denominators convergents,
$\{A_{2p+1}\}$ and $\{B_{2p+1}\}$, converge to finite limits $F_1$ and $G_1$, respectively.
Moreover, if $s$ is a finite limit point of the sequence $\{s_p\}$, and $\lim s_p =s$ as
$p$ tends to $\infty$ over a certain sequence $P$ of indices, then $A_{2p}$ and $B_{2p}$
converge to finite limits $F(s)$ and $G(s)$, respectively as $p$ tends to $\infty$ over $P$,
and
\[
F_1G(s)-G_1F(s)=1.
\]
If the sequence $\{s_p\}$ has two different finite limit points $s$ and $t$, then
\[
F(s)G(t)-F(t)G(s)=t-s.
\]
Finally, corresponding to values of $p$ for which $\lim s_p =\infty$, we have
\[
\lim\frac{A_{2p}}{B_{2p}}=\frac{F_1}{G_1},
\]
finite or infinite.
\end{theorem}

One similarity of this theorem to the present work is that it makes no
assumptions about the size
of the sequential closure.  It retains much of
the structure of the Stern-Stolz theorem, in  as much as it focuses
on the parity of the index of the approximants.  However, to understand
sequential closures in general,  all subsequences need to be
considered.  At any rate, Theorem \ref{swa} does not focus on the
sequential closure, but rather on loosening the $l_1$ assumption to
the subsequence odd indexed elements of the continued fraction.

\bigskip
\section{A generalization of certain Ramanujan Continued Fractions}

In this section we study the non-trivial case of Theorem \ref{T1}
in which the perturbing sequences $p_n$ and $q_n$
are geometric progressions tending to $0$. The inspiration
for this is the beautiful continued fraction \eqref{R3} of
Ramanujan. Theorem \ref{qcfmain} below covers
the loxodromic (convergent), parabolic, (convergent in this case),
as well as the elliptic (divergent) cases
simultaneously.   Another point
of this section is that it shows how Theorem \ref{T1} gives another
approach to the problem of evaluating continued fractions.  In fact it is
interesting to compare the proof
of Theorem \ref{qcfmain} to the proofs of special cases given
previously by different methods, see \cite{ABSYZ02,BMcL05,IS05}.

We first recall that  a  $_{1} \phi _{1}$ basic hypergeometric
series is defined for $|q|<1$ by {\allowdisplaybreaks
\begin{equation*} _{1} \phi _{1} (a;b;q,x)=  \sum_{n=0}^{\infty} \frac{(a;q)_{n}}
{(q;q)_{n}(b;q)_{n}}
 (-1)^{n} q^{n(n-1)/2}x^{n}.
\end{equation*}
}
For the $q$-product
notation used here, please see the introduction. Recall the
notation $\mathbb{T}_\lambda$
defined before Theorem \ref{T1}.
Now define $\mathbb{T}^*_\lambda$, the {\it parabolic unitary characteristic},
to be the map from $\widehat{\mathbb{C}}$ to $2^{\widehat{\mathbb{C}}}$,
equal to the set $\mathbb{T}_\lambda$, when
$\lambda\neq 1$, and $\mathbb{T}^*_1$ is a set
consisting of any fixed
element of $\widehat{\mathbb{C}}-\{1\}$. (The element of $\widehat{\mathbb{C}}-\{1\}$
does not matter.)
Again, as in Theorem \ref{T1}, we assume $\mathbb{T}^*_{\alpha/\beta}$
annihilates inconvenient terms when $|\alpha|\neq |\beta|$
and $\log_q(\alpha/\beta)\in\mathbb{Z}$.
Note that
in these cases, the continued fraction is convergent and thus
asymptotic to its limit, which will be given by the quotient
of the remaining terms on the right-hand side. Thus in \eqref{ram3gen2}
the restriction $\log_q(\alpha/\beta)\notin\mathbb{Z}$ is not too serious.

\begin{theorem}\label{qcfmain}
Let $q,\alpha,\beta\in\mathbb{C}$, $\alpha/\beta\in\widehat{\mathbb{C}}$ and $|q|<1$, then
{\allowdisplaybreaks
\begin{multline}\label{ram3gen}
\frac{ -\alpha \beta +xq}{\alpha +\beta +yq} \+
\frac{-\alpha \beta+xq^2}{ \alpha +\beta +yq^{2}} \+ \frac{-\alpha
\beta+xq^3}{
\alpha +\beta +yq^{3}}\+\cds =\\
\displaystyle{\frac{\left(\frac{x q}{\alpha} - \beta \right)
{}_1\phi_1 \left(\frac{- x q}{y \alpha}; \frac{\beta q}{\alpha};q,
\frac{-y q^2}{\alpha}
\right)\mathbb{T}^*_{\alpha/\beta}-\left(\frac{x q}{\beta} -
\alpha\right)\, {}_1\phi_1 \left(\frac{- x q}{y \beta}; \frac{\alpha
q}{\beta};q, \frac{-y q^2}{\beta} \right) } {\, {}_1\phi_1
\left(\frac{- x q}{y \alpha}; \frac{\beta q}{\alpha};q, \frac{-y
q}{\alpha} \right)\mathbb{T}^*_{\alpha/\beta}-\, {}_1\phi_1
\left(\frac{- x q}{y \beta}; \frac{\alpha q}{\beta};q, \frac{-y
q}{\beta} \right) } }.
\end{multline}
}
Moreover, assuming additionally that $\alpha \not = \beta$
and $\log_q(\alpha/\beta)\notin\mathbb{Z}$,
 {\allowdisplaybreaks
\begin{multline}\label{ram3gen2}
\frac{ -\alpha \beta +xq}{\alpha +\beta +yq} \+ \frac{-\alpha
\beta+xq^2}{ \alpha +\beta +yq^{2}} \+ \frac{-\alpha \beta+xq^3}{
\alpha +\beta +yq^{3}}\+\cds \+ \frac{-\alpha \beta+xq^n}{  \alpha
+\beta  +yq^{n}} \asymp\\
\displaystyle{\frac{\left(\frac{xq}{\alpha}-\beta \right)
{}_1\phi_1\left(\frac{- x q}{y \alpha}; \frac{\beta q}{\alpha};q,
\frac{-yq^2}{\alpha}
\right)(\frac{\alpha}{\beta})^{n+1}-\left(\frac{xq}{\beta}-
\alpha\right){}_1\phi_1\left(\frac{-xq}{y\beta};\frac{\alpha
q}{\beta};q,\frac{-y q^2}{\beta} \right)} {\, {}_1\phi_1
\left(\frac{- x q}{y \alpha}; \frac{\beta q}{\alpha};q, \frac{-y
q}{\alpha} \right)(\frac{\alpha}{\beta})^{n+1}-\, {}_1\phi_1
\left(\frac{- x q}{y \beta}; \frac{\alpha q}{\beta};q, \frac{-y
q}{\beta} \right) } }.
\end{multline}
}
\end{theorem}

This theorem contains an evaluation for a $q$-continued fraction under
the widest possible conditions. Notice that the inherent symmetry
between $\alpha$ and $\beta$ is explicit on both 
sides of the equations. We present Theorem \ref{qcfmain} as a model of the evaluation of $q$-continued fractions (for complex $|q|<1$).

In this theorem, we have not provided
the error term for the difference between the left and right hand sides
of \eqref{ram3gen2}. But
Theorem \ref{T1} implies that in the elliptic case (when $|\alpha|=|\beta|=1$),
this error is $O(q^n)$. In the loxodromic case ($|\alpha|\neq|\beta|$), the
error term can be computed from Corollary 11 in Chapter IV of \cite{LW92}.

Comparing this theorem to Theorem \ref{T1}, it is natural to enquire about
the values of the parameters $a$, $b$, $c$, and $d$. In fact, the proof
of Theorem \ref{qcfmain} follows the structure of Theorem \ref{T1} and
the constants are the expressed functions in the above statement.

{\bf Remark.} The need for using the notion of sequential closure instead
of the set of accumulation points can easily be seen here if one considers
the case $x=y=q=0$ and $\alpha=-\beta=1$. The approximants of the
continued fraction on the left hand side of \eqref{ram3gen} form the
sequence $\{\infty,0,\infty,0,\dots\}$. Accumulation points are defined
for {\it sets} and as a set it has no accumulation points. On the other hand,
one cannot use the {\it closure} of a sequence, since in general cases (where $q\neq 0$) the
approximants are not in the sequential closure. Note that in this trivial case,
$\alpha/\beta=-1$ so that $\mathbb{T}_{\alpha/\beta}^*=\{1,-1\}$, and the right hand side
of \eqref{ram3gen} simplifies to 
$(\mathbb{T}_{\alpha/\beta}^*+1)/(\mathbb{T}_{\alpha/\beta}^*-1)$ so that
both sides are equal as sets.

Before proceeding with the proof, we note a few simple corollaries. Theorem \ref{qcfmain}
generalizes certain well-known continued fraction evaluations.  For example,
setting $\alpha=y=0$ and $\beta=1$, dividing by $x$, changing $x$ to $x/q$,
and taking reciprocals in \eqref{ram3gen} yields the evaluation
of the important generalized Rogers-Ramanujan continued fraction:
\begin{corollary} For $x,q\in\mathbb{C}$ and $|q|<1$,
\[
1+\frac{xq}{1}\+\frac{xq^2}{1}\+\cds
=
\displaystyle{
\frac{\sum_{m\ge0}\frac{q^{m^2}x^m}{(q)_m}}{\sum_{m\ge0}\frac{q^{m^2+m}x^m}{(q)_m}}}.
\]
\end{corollary}

The next corollary generalizes Ramanujan's continued fraction \eqref{R3} with
three limits
given in the introduction.

\begin{corollary}\label{ramgentheor}
Let $\omega$ be a primitive $m$-th root of unity and
  let $\bar{ \omega} = 1/\omega$.
Let $1\leq i \leq m$. Then
{\allowdisplaybreaks
\begin{multline}\label{ramgentheoreq}
\lim_{k \to \infty}
\frac{1}{\omega + \bar{ \omega}+q}
\-
\frac{1}{\omega + \bar{ \omega}+q^2}
\-
\cds
\frac{1}{\omega + \bar{ \omega}+q^{mk+i}}\\
= \displaystyle{\frac{ \omega^{1-i} \,_1 \phi_1 \left(0;
 q\omega^2;q, -q^2\omega \right)- \omega^{i-1} \,_1 \phi_1 \left(0;
q/ \omega^2 ;q, -q^2/\omega \right) } { \omega^{-i} \,_1 \phi_1
\left(0;
 q\omega^2;q, -q\omega \right)- \omega^{i} \,_1 \phi_1 \left(0;
q/ \omega^2 ;q, -q/\omega \right) } }.
\end{multline}
}
\end{corollary}
\begin{proof}
This is immediate from \eqref{ram3gen2}, upon setting $x=0$, $y=1$,
$\alpha = \omega$, $\beta = \omega^{-1}$, $n=mk+i$, then noting that
$\omega^{mk}=1$.
\end{proof}

This result in its present form first appeared in \cite{IS05}. The authors
found it independently and gave a different proof
in \cite{BMcL05}.

\begin{corollary}\label{32cor}
\begin{equation*}
\frac 32 - \frac 1{3/2}\- \frac 1{3/2}\- \frac 1{3/2}\-\cds =\mathbb{R}.
\end{equation*}
In fact, the limiting probability density function of the approximants of this continued fraction
is given by
\[
p(x) =\frac{\sqrt{7}}{2 \pi(2 x^2-3 x+2)}.
\]
\end{corollary}
\begin{proof} In Theorem \ref{qcfmain} take $\alpha=3/4+i\sqrt{7}/4$,
$\beta=3/4-i\sqrt{7}/4$, $q=x=y=0$. The limiting probability density function follows
from Corollary \ref{c3}.
\end{proof}

\begin{figure}[htbp]
\centering \epsfig{file=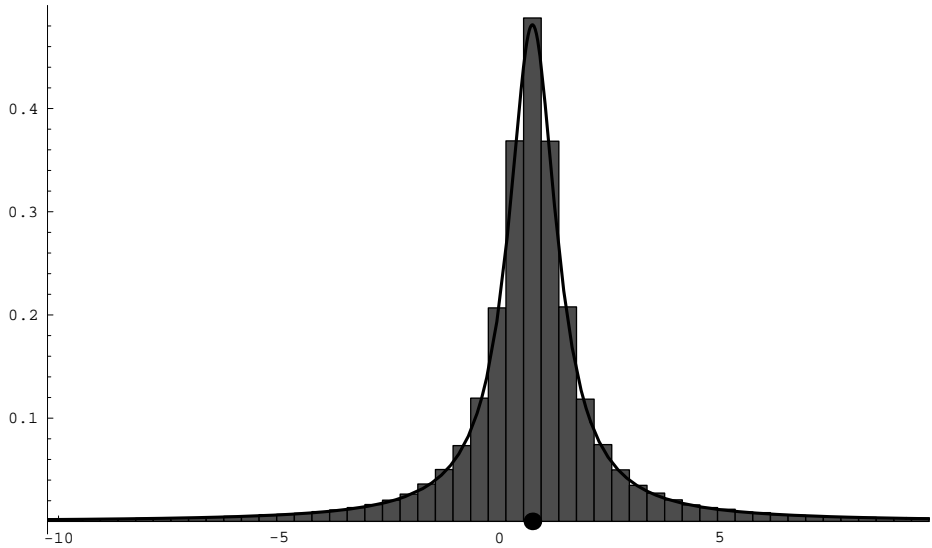, width= 350pt}
\caption{The distribution of the first 3000 approximants of
$3/2+\mathop{\pmb{K}}_{n=1}^{\infty}\frac{-1}{3/2}$, with the point $x=3/4$ of
predicted highest concentration and the limiting probability density
function $p(x)=\sqrt{7}/(2 \pi(2 x^2-3 x+2))$. } \label{fig4a}
\end{figure}

Thus Theorem \ref{qcfmain} unifies  the celebrated
Rogers-Ramanujan continued fraction,  Ramanujan's
continued fraction with three limits, and \eqref{32R}, the continued
fraction for $\mathbb{R}$ given in the introduction;  it gives identities that can
have any rank in $\mathbb{Z}^+\cup \{\mathfrak{c}\}$.

It is interesting to consider that although identities such as that
in Corollary \ref{32cor} may appear useless as they display the
oscillatory divergence of the continued fraction, the divergence is
actually quite well controlled in the sense that there is an
explicit probability density function for the approximants.
Later Corollary \ref{c3} will be used to compute the
point of highest concentration as the average of two convergent
variants of the original continued fraction. Thus these continued
fractions come with a kind of built-in summability. Of course, if
one is interested in computing the sequential closure circle, the
continued fraction converge to it rather rapidly, as was remarked in
the first paragraph following Theorem \ref{qcfmain}. See also Figure
\ref{fig3} below.

Here is the proof of Theorem \ref{qcfmain}. Following the proof
other special cases are studied.

\begin{proof}
This theorem is a simple consequence of Theorem \ref{T1} and work
from our paper \cite{BMcLW06}.
First consider the case $|\alpha|=|\beta|=1$.  Let $P_n$ and $Q_n$
denote the $n$th numerator and denominator convergents of \eqref{ram3gen2}.
In \cite{BMcLW06} we worked with the related continued fraction
 {\allowdisplaybreaks
\begin{equation}\label{cfeq}
\frac{1}{1} \+ \frac{-\alpha \beta+xq }{\alpha +\beta +yq} \+
\frac{-\alpha \beta+xq^2}{  \alpha +\beta +yq^{2}} \+ \frac{-\alpha
\beta+xq^3}{  \alpha +\beta  +yq^{3}}\+ \cds.
\end{equation}
}
(Note that in \cite{BMcLW06} $a$, $b$, $c$, and $d$ were written
for $\alpha$, $\beta$, $x$, and $y$, respectively.)
Let $A_{n}$ and $B_{n}$ denote the $n$-th numerator convergent
and $n$-th denominator convergent, respectively, of \eqref{cfeq}.
Then it is immediate that $P_n=B_{n+1}-A_{n+1}$ and
$Q_n=A_{n+1}$. Observe that since both of these continued fractions
are symmetric in $\alpha$ and $\beta$, the explicit
expressions for $A_n$ and $B_n$ ((2.5-6) from \cite{BMcLW06})
are also
valid with $\alpha$ and $\beta$ interchanged. All that is needed
is to compute the limits \eqref{abcddef} from Theorem \ref{T1}.
Now employing the convergents of \eqref{ram3gen2} instead
of those for \eqref{cfeq} puts equations (2.15) and (2.16) from
\cite{BMcLW06} into the forms

\begin{equation}\label{soon}
\lim_{N\to\infty}\frac{P_N}{\beta^{N}}
=
\frac{(xq/\beta-\alpha)}{1-\alpha/\beta}\sum_{n=1}^\infty
\frac{(-xq/y\beta)_n}{(\alpha q/\beta)_n(q)_n}(-1)^nq^{n(n-1)/2}\left(\frac{-yq^2}{\beta}\right)^n,
\end{equation}
and
\begin{equation}\label{sooner}
\lim_{N\to\infty}\frac{Q_N}{\beta^{N}}
=\frac{1}{1-\alpha/\beta}
\sum_{n=1}^\infty
\frac{(-xq/y\beta)_n}{(\alpha q/\beta)_n(q)_n}(-1)^nq^{n(n-1)/2}\left(\frac{-yq}{\beta}\right)^n.
\end{equation}
Thus, for example, using \eqref{abcddef},
\begin{equation*}
\begin{split}
b&=-\lim_{N\to\infty}\beta^{-N}(P_N-\alpha P_{N-1})\\
&=-\left(\lim_{N\to\infty}\beta^{-N}P_N-(\alpha/\beta)\lim_{N\to\infty}\beta^{-(N-1)}P_{N-1}\right)\\
&=-\left(1-\alpha/\beta\right)\lim_{N\to\infty}\beta^{-N}P_N\\
&=-\left(\frac{x q}{\beta} -
\alpha\right)\,_1 \phi_1 \left(\frac{- x q}{y \beta}; \frac{\alpha
q}{\beta};q, \frac{-y q^2}{\beta} \right).
\end{split}
\end{equation*}
Here the last three equalities followed from \eqref{soon}.

The identification of the other constants in $h$ is similar, except that
one must interchange the role of $\alpha$ and $\beta$ in \eqref{soon}
and \eqref{sooner} when calculating $a$ and $c$.

The case $|\alpha|=|\beta|\neq 0$ follows by taking the
equivalence transformation as in the proof of Theorem \ref{T1}.
Note that the parameters $\alpha$ and $\beta$ in the $\phi$
function are not restricted to this domain. For brevity in
the rest of the proof, we employ the variables $a$, $b$,
$c$, and $d$, with the understanding that they represent
the $\phi$ functions with the above arguments, analytically
continued to their widest domain. (It is easy to check that
the $\phi$ function is meromorphic in its variables in the
complex plane.)

Now assume that $|\alpha|\neq|\beta|$. First note that the
difference equation
\begin{equation}\label{qde}
Y_n=(1+\lambda-zq^{n})Y_{n+1}+(-\lambda+szq^n)Y_{n+2}
\end{equation}
has a solution $Y_n=\,_1\phi_1(s;\lambda q;q,zq^n)$. (This can be checked simply by equating
coefficients.) By Auric's theorem,
see Corollary 11, Chapter IV of \cite{LW92}, this solution of \eqref{qde}
is minimal if $|\lambda|<1$ or $\lambda=1$, and thus for $|\lambda|<1$
or $\lambda=1$,
\[
\frac{\,_1\phi_1(s;\lambda q;q,z)}{\,_1\phi_1(s;\lambda q;q,zq)}=1+\lambda-z+
\frac{-\lambda+sz}{1+\lambda-zq}\+\frac{-\lambda+szq}{1+\lambda-zq^2}\+\cds.
\]
Putting $s=-\beta^{-1}xy^{-1}q$, $\lambda=\alpha/\beta$, and $z=-\beta^{-1}yq$,
taking reciprocals, multiplying both sides by $-\alpha+xq/\beta$
and applying a simple equivalence transformation to the continued fraction,
yields that for $|\alpha|<|\beta|$ or $\alpha=\beta\neq 0$,
 {\allowdisplaybreaks
\begin{equation}\label{twice}
\displaystyle{
\frac{b}{d}
 }
=
\frac{ -\alpha \beta +xq}{\alpha +\beta +yq} \+ \frac{-\alpha
\beta+xq^2}{ \alpha +\beta +yq^{2}} \+ \frac{-\alpha \beta+xq^3}{
\alpha +\beta +yq^{3}}\+\cds .
\end{equation}
}
For $|\alpha|>|\beta|$, symmetry gives that
 {\allowdisplaybreaks
\begin{equation}\label{thrice}
\displaystyle{
\frac ac
 }
=
\frac{ -\alpha \beta +xq}{\alpha +\beta +yq} \+ \frac{-\alpha
\beta+xq^2}{ \alpha +\beta +yq^{2}} \+ \frac{-\alpha \beta+xq^3}{
\alpha +\beta +yq^{3}}\+\cds .
\end{equation}
}
The conclusion follows by noting that for $|\alpha|<|\beta|$ as $n\to\infty$,
 {\allowdisplaybreaks
\begin{equation*}
\frac{a\lambda^{n+1}+b}{c\lambda^{n+1}+d}\asymp\frac bd,
\end{equation*}}
while for $|\alpha|>|\beta|$,
 {\allowdisplaybreaks
\begin{equation*}
\frac{a\lambda^{n+1}+b}{c\lambda^{n+1}+d}\asymp\frac ac.\end{equation*}
}
\end{proof}

{\bf Remark.} We could have simply used \eqref{soon} and
\eqref{sooner} to complete the $|\alpha|\neq|\beta|$ part of the
proof, but the approach via Auric's theorem seemed preferable as it
also yields the evaluation of the continued fraction in the
parabolic case $\alpha=\beta\neq 0$.

Consider the special case of the continued fraction in the theorem in
which $x=0$ and $y=1$. Then
\begin{equation}\label{hzeqa}
h(z) =\displaystyle{\frac{- \beta \,_1 \phi_1 \left(0; \frac{\beta
q}{\alpha};q, \frac{- q^2}{\alpha} \right)z+ \alpha \,_1 \phi_1
\left(0; \frac{\alpha q}{\beta};q, \frac{- q^2}{\beta} \right) }
{\,_1 \phi_1 \left(0; \frac{\beta q}{\alpha};q, \frac{- q}{\alpha}
\right)z-\,_1 \phi_1 \left(0; \frac{\alpha q}{\beta};q, \frac{-
q}{\beta} \right) } },
\end{equation}
and thus that the sequential closure of the continued fraction
\begin{equation*}
G(\alpha, \beta, q):=\frac{ -\alpha \beta }{\alpha +\beta +q} \-
\frac{\alpha \beta}{  \alpha +\beta +q^{2}} \- \frac{\alpha \beta}{
\alpha +\beta  +q^{3}} \cds
\end{equation*}
is on the circle $f(\mathbb{T})$, where $f$ is defined by
\begin{equation*}
f(z)= \displaystyle{\frac{- \beta \,_1 \phi_1 \left(0; \frac{\beta
q}{\alpha};q, \frac{- q^2}{\alpha} \right)z+ \alpha \,_1 \phi_1
\left(0; \frac{\alpha q}{\beta};q, \frac{- q^2}{\beta} \right) }
{\,_1 \phi_1 \left(0; \frac{\beta q}{\alpha};q, \frac{- q}{\alpha}
\right)z-\,_1 \phi_1 \left(0; \frac{\alpha q}{\beta};q, \frac{-
q}{\beta} \right) }   }.
\end{equation*}

\begin{figure}[htbp]
\centering \epsfig{file=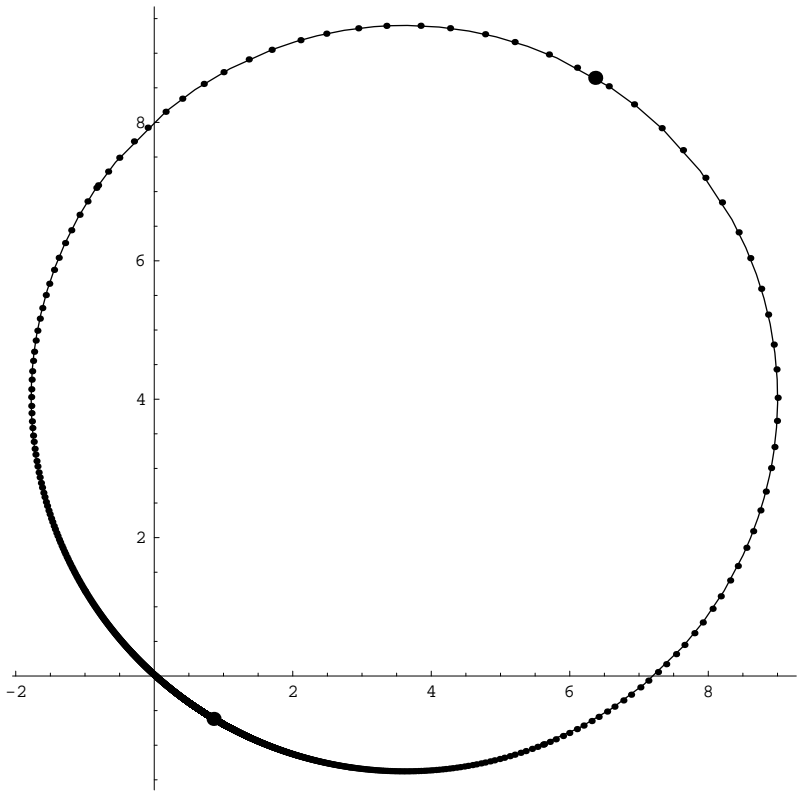, width= 300pt} \caption{The
convergence of $G(\exp(\imath \sqrt{7}),\exp(\imath  \sqrt{5}),
0.1)$ } \label{fig3}
\end{figure}

\bigskip

Figure \ref{fig3} shows the first 3100 approximants of
$G(\exp(\imath \sqrt{7}),\exp(\imath  \sqrt{5}), 0.1)$ and the
corresponding circle $f(\mathbb{T})$ predicted by the theory. The
larger dots show the points, again predicted by the theory, of
highest and lowest concentration of approximants. Note the error,
$\varepsilon_n=O(10^{-n})$ and experimentally, $\min_{z \in
\mathbb{T}}|A_{n}/B_{n} - f(z)| \approx 10^{-n}$ in agreement with
the theory.

Figure \ref{fig4} shows the first 2700 approximants (with the first
55 omitted, since they lay relatively distant from the circle of
convergence) of\\ $G(\exp(\imath \sqrt{7}),\exp(\imath (\sqrt{7}+2
\pi/11)) , 0.99\exp(\imath \sqrt{17}))$ and its convergence to the
eleven limit points $f(2 k \pi/11)$, $0 \leq k \leq 10$, where
$f(z)$ is the associated linear fractional transformation, together
with the circle $f(\mathbb{T})$. The value $q=0.99\exp(\imath
\sqrt{17}))$ was chosen to be close to 1 in absolute value, with the
aim of slowing down the convergence so as to make the behavior more
visible.

{\allowdisplaybreaks
\begin{figure}[htbp]
\centering \epsfig{file=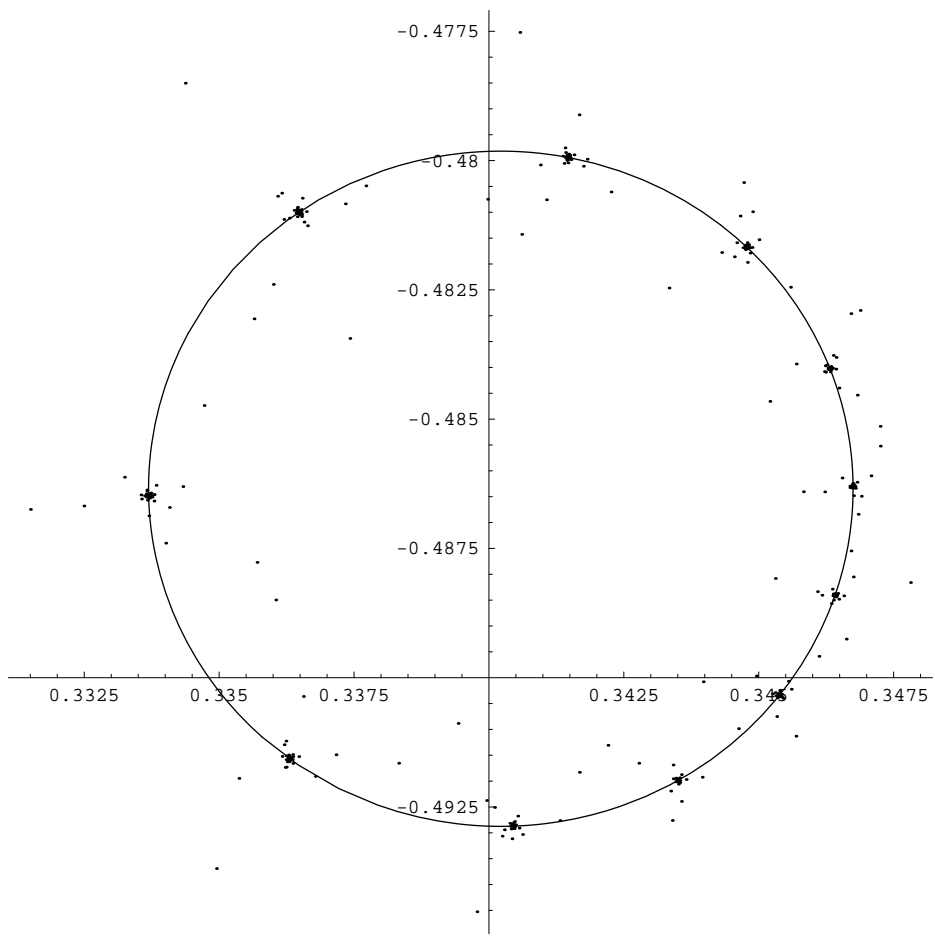, width= 350pt}
\caption{The convergence of $G(\exp(\imath \sqrt{7}),\exp(\imath
(\sqrt{7}+2 \pi/11)) , 0.99\exp(\imath \sqrt{17}))$ } \label{fig4}
\end{figure}
}

The Bauer-Muir transformation (Proposition \ref{BM}), can
applied to Theorem \ref{qcfmain}
to obtain convergent
continued fractions.  The results are contained in the following corollary.
The continued fraction from Theorem \ref{qcfmain}
will be denoted as follows.
\[
K(\alpha,\beta,x,y):=
\frac{ -\alpha \beta +xq}{\alpha +\beta +yq} \+
\frac{-\alpha \beta+xq^2}{ \alpha +\beta +yq^{2}} \+ \frac{-\alpha
\beta+xq^3}{
\alpha +\beta +yq^{3}}\+\cds.
\]

\begin{corollary}\label{rBM}
Let $|q|<1$ and let $|\alpha|=|\beta|\neq 0$ be distinct complex
numbers such that $\arg(\alpha/\beta)$ is not a rational multiple of
$\pi$. When the approximants of $K(\alpha,\beta,x,y)$ are real, so
that the continued fraction  is dense on $\mathbb{R}$, its point of
highest concentration is given by {\allowdisplaybreaks
\begin{multline*} x_0= \\\frac 12\left(-\beta+ \frac{\beta yq+xq}{\alpha
+yq+K(\alpha,\beta q,xq,yq)}-\alpha+ \frac{ \alpha yq+xq}{\beta
+yq+K(\alpha q,\beta,xq,yq)}\right).
\end{multline*}
}
 In fact, the approximants of $K(\alpha,\beta,x,y)$ have the Cauchy
distribution
\[
p(x)=\frac{\delta}{\pi[(x-x_0)^2+\delta^2]},
\]
with scale factor
\[
\delta=
\frac 1{2i}\left(-\beta+
\frac{ \beta yq+xq}{\alpha  +yq+K(\alpha,\beta q,xq,yq)}
 +\alpha-
\frac{ \alpha yq+xq}{\beta  +yq+K(\alpha q,\beta ,xq,yq)}    \right).
\]

Moreover, assuming only the first sentence,
\begin{equation}\label{hinq}
 {\allowdisplaybreaks
h(\infty)=
\displaystyle{
\frac{\left(\frac{x q}{\alpha}-\beta \right)\,
{}_1\phi_1 \left(\frac{- x q}{y \alpha}; \frac{\beta q}{\alpha};q,
\frac{-y q^2}{\alpha} \right)}{{}_1\phi_1
\left(\frac{- x q}{y \alpha}; \frac{\beta q}{\alpha};q, \frac{-y
q}{\alpha} \right)}
 }
= -\beta+
\frac{ \beta yq+xq}{\alpha  +yq+K(\alpha,\beta q,xq,yq)} ,}
\end{equation}
and,
\begin{equation}\label{hzrq}
 {\allowdisplaybreaks
h(0)=
\displaystyle{
\frac{\left(\frac{x q}{\beta}-\alpha \right)\,
{}_1\phi_1 \left(\frac{- x q}{y \beta}; \frac{\alpha q}{\beta};q,
\frac{-y q^2}{\beta} \right)}{{}_1\phi_1
\left(\frac{- x q}{y \beta}; \frac{\alpha q}{\beta};q, \frac{-y
q}{\beta} \right)}
 }
= -\alpha+
\frac{ \alpha yq+xq}{\beta  +yq+K(\alpha q,\beta ,xq,yq)} .}
\end{equation}
\end{corollary}

This corollary exhibits the interesting phenomenon of a continued
fraction that diverges (and is dense in $\mathbb{R}$), yet the
statistical distribution is given as a rational function of the same
continued fraction evaluated at points where it converges.

\begin{proof}
By Corollary \ref{modlimits},  it follows that
$h(\infty)=\lim_{n\to\infty}f_n(-\beta)$ and $h(0)=\lim_{n\to\infty}f_n(-\alpha)$,
where $f_n(w)$ denotes the $n$th modified approximant of $K(\alpha,\beta,x,y)$.
Apply Proposition \ref{BM} to the continued fraction $K(\alpha,\beta,x,y)$
with the modifications $w_n$ taken to be the constant sequences
$\{-\beta, -\beta,\dots\}$ and $\{-\alpha, -\alpha,\dots\}$, respectively,
to obtain the convergent continued fractions for $h(\infty)$ and $h(0)$.
The series expressions for $h(\infty)$ and $h(0)$ follow directly from Theorem \ref{qcfmain}.
The final conclusion is then an application of Corollary \ref{c3}.  The conditions
$x\neq -\beta y$ and $x\neq -\alpha y$ are needed to ensure that
the Bauer-Muir transformations exist ($\lambda_n\neq 0$ in Proposition \ref{BM}).
Observe that these conditions may be dropped in the corollary
by the following well-known
version of the $q$-binomial theorem \cite{GR}:
\[
\sum_{n=0}^\infty \frac{q^{(n-1)n/2}z^n}{(q)_n}=(-z)_\infty.
\]
\end{proof}

Notice that the two continued fraction to series identities in this
corollary are equivalent by the symmetry in $\alpha$ and $\beta$. Also,
the series equals continued fraction identities follow from Theorem \ref{qcfmain}
by setting $\beta$ to $\beta q$ and then assuming $|\alpha|=|\beta|\neq 0$,
and manipulation of the continued fraction and series. (A contiguous relation
needs to be employed to simplify the series in this case.)

In some cases the infinite series in the corollary above can be
expressed as infinite products.

\begin{corollary} Let $|q|<1$. Then
\begin{equation}\label{BMRcf3a}
1 -\frac{q}{1+q} \+
 \mathop{\pmb{K}}_{n=2}^{\infty}\frac{q^2}
   { 1 -q^2+q^{2n-1} }
   = \frac{(q;q^5)_{\infty}(q^4;q^5)_{\infty}
   }
   {(q^2;q^5)_{\infty}(q^3;q^5)_{\infty}
   }.
\end{equation}
\end{corollary}

\begin{proof}
In \eqref{hinq}, set $x=0$, $y=1$, replace $q$ by $q^2$, set $\beta = -q$ and
$\alpha =q$ and simplify the resulting continued fraction by
applying a sequence of equivalence transformations.

For the right side we use two identities due to Rogers
\cite{Rogers} (see also \cite{S52} and \cite{S03}, identities
\textbf{A.16} and \textbf{A.20}):
\begin{align*}
\sum_{n=0}^{\infty}\frac{  q^{n(n+2)}}{(q^4;q^4)_{n}}
&=\frac{1}{(q^2;q^5)_{\infty}(q^3;q^5)_{\infty}(-q^2;q^2)_{\infty}},\\
\sum_{n=0}^{\infty}\frac{q^{n^2}}{(q^4;q^4)_{n}}
&=\frac{1}{(q;q^5)_{\infty}(q^4;q^5)_{\infty}(-q^2;q^2)_{\infty}}.
\end{align*}
Finally, cancel a factor of $q$ on each side
\end{proof}
The continued fraction above is clearly a transformed
version of the Rogers-Ramanujan continued fraction since they converge
to the same limit.

\section{Applications to $(r,s)$-matrix continued fractions}

Levrie and Bultheel   defined a generalization of continued
fractions called $(r,s)$-matrix continued fractions \cite{LB}. This
generalization unifies a number of generalizations of continued
fractions including ``generalized (vector valued) continued
fractions" and continued fractions obtained by composing together
higher dimensional M\"obius transformations, see \cite{LW92}
and \cite{Bear1} for definitions.

Here we show that our results apply to limit periodic
$(r,s)$-matrix continued fractions with eigenvalues of equal
magnitude yielding estimates for the asymptotics of their
approximants.

For consistency we closely follow the notation used in \cite{LB} to define $(r,s)$-matrix
continued fractions. Let $\mathfrak{M}_{s,r}(\mathbb{C})$ denote the set of $s\times r$ matrices
over the complex numbers.  Let $\theta_k$ be a sequence of $n\times n$ matrices over
$\mathbb{C}$. Assume that $r+s=n$.  A $(r,s)$-matrix continued fraction is associated
with a recurrence system of the form $Y_k=Y_{k-1}\theta_k$. The continued
fraction is defined by its sequence of approximants. These are sequences of $s\times r$
matrices defined in the following manner.

Define the function $f: D\in \mathfrak{M}_n(\mathbb{C})\to
\mathfrak{M}_{s,r}(\mathbb{C})$ by
\begin{equation}\label{fmateq} f(D)=
B^{-1}A,
\end{equation}
where $B$ is the $s\times s$ submatrix
consisting of the last $s$ elements from both the rows and columns
of $D$, and $A$ is the $s\times r$ submatrix consisting of the
first $r$ elements from the last $s$ rows of $D$.

Then the $k$-th approximant of the $(r,s)$-matrix continued
fraction associated with the sequence $\theta_k$ is defined to be
\begin{equation}\label{skeq}
s_k:=f(\theta_k\theta_{k-1}\cdots\theta_2\theta_1).
\end{equation}

To apply Theorem \ref{gent} to this situation, we endow
$\mathfrak{M}_{s\times r}(\mathbb{C})$ with a metric by letting the distance
function for two such matrices be the maximum absolute value of
the respective differences of corresponding pairs of elements.
Then when $f$ is continuous, our theorem can be
applied. (Note that $f$ will be continuous provided that it
exists.)

Let $\lim_{k\to\infty}\theta_k=\theta$, for some $\theta\in \mathfrak{M}_n(\mathbb{C})$. Then the
recurrence system is said to be of Poincar{\'e} type and the
$(r,s)$-matrix continued fraction is called limit periodic.

After this definition Theorem \ref{gent} can be applied and the
following theorem results.

\begin{theorem}\label{rs1}
Suppose that the condition $\sum_{k\ge
1}||\theta_k-\theta||<\infty$ holds, that the matrix $\theta$ is
diagonalizable, and that the eigenvalues of $\,\theta$ are all of
magnitude $1$. Then the $k$th approximant $s_k$ has the asymptotic
formula
\begin{equation}\label{sksim}
s_k\asymp f(\theta^k F),
\end{equation}
where $F$ is the matrix defined by the convergent product
\[
F:=\lim_{k\to\infty}\theta^{-k}\theta_k\theta_{k-1}\cdots\theta_2\theta_1.
\]
\end{theorem}

Note that because of the way that $(r,s)$-matrix continued fractions are defined,
we have taken products in the reverse order than the rest of the paper.

As a consequence of this asymptotic, the sequential closure can be determined from
\[
\copyright(s_k)=\copyright(f(\theta^k F)).
\]

In one  general case, detailed in the following theorem, we
actually get a convergence theorem.
\begin{theorem}\label{rs2}
Let $\theta_k$ be a sequence of $n\times n$ matrices over
$\mathbb{C}$ satisfying
\[
\sum_{k\ge 1}||\theta_k-\theta||<\infty, \] where $\theta$ is a
diagonal (or antidiagonal) matrix with all diagonal (or antidiagonal)
elements of absolute value 1.
Let $r$ and $s$ be positive integers with $r+s=n$.

Then the matrix
\[
F:=\lim_{k\to\infty}\theta^{-k}\theta_k\theta_{k-1}\cdots\theta_2\theta_1
\]
exists. Suppose further that the bottom right $s \times s$
submatrix of $F$ is nonsingular. Then the $(r,s)$-matrix continued
fraction defined by equation \eqref{skeq} converges to $f(F)$. If $\theta$
is antidiagonal, then the even approximants of the $(r,s)$-matrix continued
fraction defined by equation \eqref{skeq} tend to $f(F)$, while its odd approximants
tend to $f(AF)$,
where $A$ is the antidiagonal matrix with $1$s along its antidiagonal.
\end{theorem}
\begin{proof}
The matrix $F$ exists by Theorem \ref{gent} (or more precisely,
the ``transposed" version of Theorem \ref{gent}). Let
\[
\theta = \textnormal{diag} (\lambda_1, \dots, \lambda_n).
\]
By \eqref{sksim},
\begin{align*}
s_k &\asymp f(\theta^k F)\\ &=\left(\left (
\begin{matrix}
&\lambda^k_{n-s+1} &\dots &0\\
&\vdots& \ddots &\vdots\\
&0 &\dots & \lambda^k_n
\end{matrix}
\right) \left (
\begin{matrix}
&F_{n-s+1,n-s+1}&\dots &F_{n-s+1,n}\\
&\vdots &\ddots &\vdots\\
&F_{n,n-s+1} &\dots  & F_{n,n}
\end{matrix}
\right)\right)^{-1}\\
&\phantom{asdadasd}\times \left (
\begin{matrix}
&\lambda^k_{n-s+1} &\dots &0\\
&\vdots& \ddots &\vdots\\
&0 &\dots & \lambda^k_n
\end{matrix}
\right) \left (
\begin{matrix}
&F_{n-s+1,1}&\dots &F_{n-s+1,r}\\
&\vdots &\ddots &\vdots\\
&F_{n,1} &\dots  & F_{n,r}
\end{matrix}
\right)\\
&=\left (
\begin{matrix}
&F_{n-s+1,n-s+1}&\dots &F_{n-s+1,n}\\
&\vdots &\ddots &\vdots\\
&F_{n,n-s+1} &\dots  & F_{n,n}
\end{matrix}
\right)^{-1} \left (
\begin{matrix}
&F_{n-s+1,1}&\dots &F_{n-s+1,r}\\
&\vdots &\ddots &\vdots\\
&F_{n,1} &\dots  & F_{n,r}
\end{matrix}
\right)\\
&=f(F).
\end{align*}
Thus $s_k$ converges to the final matrix product above.

For the case where $\theta$ is an antidiagonal matrix, $\theta^{2k}$
is a diagonal matrix and the proof for the even approximants is
virtually the same as for the case where $\theta$ is a diagonal
matrix. If $\theta$ is an antidiagonal matrix, $\theta^{2k+1}$ is
also an antidiagonal matrix. Once again by \eqref{sksim},
{\allowdisplaybreaks
\begin{align*}
&s_{2k+1} \asymp f(\theta^{2k+1} F)\\ &=\left(\left (
\begin{matrix}
&0 &\dots &(\theta^{2k+1})_{n-s+1,s}\\
&\vdots& \iddots &\vdots\\
&(\theta^{2k+1})_{n,1} &\dots & 0
\end{matrix}
\right) \left (
\begin{matrix}
&F_{1,n-s+1}&\dots &F_{1,n}\\
&\vdots &\ddots &\vdots\\
&F_{s,n-s+1} &\dots  & F_{s,n}
\end{matrix}
\right)\right)^{-1}\\
&\phantom{asdadasd}\times \left (
\begin{matrix}
&0 &\dots &(\theta^{2k+1})_{n-s+1,s}\\
&\vdots& \iddots &\vdots\\
&(\theta^{2k+1})_{n,1} &\dots & 0
\end{matrix}
\right) \left (
\begin{matrix}
&F_{1,1}&\dots &F_{1,r}\\
&\vdots &\ddots &\vdots\\
&F_{s,1} &\dots  & F_{s,r}
\end{matrix}
\right)\\
&=\left(\left (
\begin{matrix}
&0 &\dots &1\\
&\vdots& \iddots &\vdots\\
&1 &\dots & 0
\end{matrix}
\right) \left (
\begin{matrix}
&F_{1,n-s+1}&\dots &F_{1,n}\\
&\vdots &\ddots &\vdots\\
&F_{s,n-s+1} &\dots  & F_{s,n}
\end{matrix}
\right)\right)^{-1}\\
&\phantom{asdsadasdadasdasdasdaddadasd}\times \left (
\begin{matrix}
&0 &\dots &1\\
&\vdots& \iddots &\vdots\\
&1 &\dots & 0
\end{matrix}
\right) \left (
\begin{matrix}
&F_{1,1}&\dots &F_{1,r}\\
&\vdots &\ddots &\vdots\\
&F_{s,1} &\dots  & F_{s,r}
\end{matrix}
\right)
\\
&=f(AF),
\end{align*}
}where $A$ is the antidiagonal matrix with $1$'s along the
antidiagonal. Thus $s_{2k+1}$ converges to the final matrix product
above.
\end{proof}

Consider now the $n=2$ antidiagonal case of Theorem \ref{rs2}. The matrix $\theta$
then has the form

\[
\theta=
\left(
\begin{matrix}
&0 &1\\
&1 &0
\end{matrix}
\right).
\]
Choose $\theta_k$ to have the form

\[
\theta_k=
\left(
\begin{matrix}
&b_k &1\\
&1+a_k &0
\end{matrix}
\right).
\]

Using the correspondence between matrices and continued fractions \eqref{corr}, we at once obtain
the following corollary, first given in \cite{BMcL05}.

\begin{corollary}\label{SSgen}
Let the sequences $\{a_n\}$ and $\{b_{n}\}$ satisfy $a_n\ne -1$ for $n\ge 1$,  $\sum |a_{n}| <\infty$ and
$\sum |b_{n}| <\infty$.
Then
\[
b_{0}+\mathop{\pmb{K}}_{n=1}^{\infty}\frac{1+a_n}{b_{n}}
\]
diverges. In fact, for $p=0,1$,
\begin{align*}
&\lim_{n \to \infty}P_{2n+p}=A_{p} \not = \infty,& &\lim_{n \to \infty}Q_{2n+p}=B_{p} \not = \infty,&
\end{align*}
and
\[
A_{1}B_{0}-A_{0}B_{1} = \prod_{n=1}^{\infty}(1+a_n).
\]
\end{corollary}

In fact, Corollary \ref{SSgen} is also the
$\alpha=1$,
$\beta = -1$ (so $m=2$), $q_{n}=a_n$ and $p_{n} = b_{n}$
 case of Corollary \ref{c1}.  When $a_n=0$, this corollary
reduces to the famous Stern-Stolz theorem discussed in the introduction.

One of the main results of the paper \cite{BMcL05} was Corollary \ref{c1}, which
we applied to obtain an infinite sequence of theorems, similar to the Stern-Stolz
 theorem, but with continued fractions of different ranks.
Notice that Theorem \ref{rs2} provides yet another family of generalizations.

It is interesting to compare Corollary \ref{SSgen} with the ``The General
Stern-Stolz Theorem'' from \cite{Bear1} in the case of continued fractions.
The corollary for the case of complex continued fractions is:

\begin{corollary}\label{SSbiggen}[Corollary 7.5 of \cite{Bear1}]
If $\sum_n |1-|a_n||$ and $\sum_n|b_n|$ converge, then $\mathop{\pmb{K}}(a_n|b_n)$ is
strongly divergent.
\end{corollary}

The first condition in this result is weaker than analogous
condition in Corollary \ref{SSgen} above. But it should be remarked
that Theorem \ref{Stern-Stolz1}, Corollary \ref{SSgen}, and
Corollary \ref{SSbiggen} are, in fact, equivalent;  the two
corollaries follow from Theorem \ref{Stern-Stolz1} by an equivalence
transformation (and a little analysis). Next, the condition on the
partial numerators in Corollary \ref{SSbiggen} encodes the
information that the matrices representing the continued fraction
are a perturbation of unitary matrices. We could have obtained the
same result by using Theorem \ref{DMunitconj}, however in this
situation one does not obtain as detailed information about the
limits of the convergents. In particular, Corollary \ref{SSgen} also
proves the convergence of the subsequences of convergents $\{P_n\}$
and $\{Q_n\}$ of equal parity. Corollary \ref{SSbiggen} does not
furnish this part of the conclusion. On the other hand, it does
prove strong divergence, defined in section 2. Indeed, the continued
fraction in Corollary \ref{SSbiggen} is not necessarily limit
periodic.

\section{Poincar{\'e} type recurrence relations with characteristic roots on the unit circle}

Let the complex sequence $\{x_{n}\}_{n \geq 0}$ have the initial values
$x_{0}$, $\dots$, $x_{p-1}$ and subsequently be defined by
\begin{equation}\label{pereq}
x_{n+p}=\sum_{r=0}^{p-1}a_{n,r}x_{n+r},
\end{equation}
for $n \geq 0$. It is assumed that for $n$ sufficiently large $a_{n,0}a_{n,p-1}\neq 0$.
Suppose also that there are numbers $a_{0},\dots, a_{p-1}$ such that
\begin{align}
\label{asumineq}
&\lim_{n \to \infty}a_{n,r}=a_{r},& &0 \leq r \leq p-1.&
\end{align}

A recurrence of the form (\ref{pereq}) satisfying the condition
(\ref{asumineq}) is called a Poincar{\'e}-type recurrence. Such
recurrences were initially studied by\\ Poincar{\'e} and later
Perron who proved the Poincar\'e-Perron theorem \cite{koomR,Poi}:

{\it If the roots of the characteristic equation
\begin{equation}\label{chareq}
t^{p}-a_{p-1}t^{p-1}-a_{p-2}t^{p-2}- \dots -a_{0}=0
\end{equation}
have distinct norms, then $\lim_{n\to\infty}x_{n+1}/x_n=\alpha$,
where $\alpha$ is a root of \eqref{chareq}. Moreover, for each
root $\alpha$ of \eqref{chareq}, there exists
a solution of \eqref{pereq} with term ratio tending to $\alpha$.}

 Because the roots are also the
eigenvalues of the associated companion matrix, they  are also
referred to as the eigenvalues of (\ref{pereq}).  This result was
improved by O. Perron, who obtained a number of theorems about the
limiting asymptotics of such recurrence sequences. Perron
\cite{Perron1} made a significant advance in 1921 when he proved
the following theorem which treated cases of
eigenvalues which repeat or are of equal norm.

\begin{proposition}\label{tp1}
Let the sequence $\{x_{n}\}_{n \geq 0}$ be defined by initial values
$x_{0}$, $\dots$, $x_{p-1}$ and by (\ref{pereq})
for $n \geq 0$. Suppose also that there are numbers $a_{0},\dots, a_{p-1}$ satisfying (\ref{asumineq}).
Let $q_{1},\,q_{2}, \dots q_{\sigma}$ be the distinct moduli of the roots of the characteristic
equation (\ref{chareq})
and let $l_{\lambda}$ be the number of roots whose modulus is $q_{\lambda}$, multiple roots counted
according to multiplicity, so that
\[
l_{1}+l_{2}+\dots l_{\sigma}=p.
\]
Then, provided $a_{n,0}$ be different from zero for $n \geq 0$, the difference equation
\eqref{pereq} has a fundamental system of solutions,
which fall into $\sigma$ classes, such that, for the solutions of the $\lambda$-th class
and their linear combinations,
\[
\limsup_{n \to \infty}  \sqrt[n]{|x_{n}|}=q_{\lambda}.
\]
The number of solutions of the $\lambda$-th class is $l_{\lambda}$.
\end{proposition}

Thus when all of the characteristic roots have norm $1$, this theorem
gives that
\[
\limsup_{n  \to \infty}  \sqrt[n]{|x_{n}|}= 1.
\]

Another related paper is \cite{KL} where the authors study products
of matrices and give a sufficient condition for their boundedness. This
is then used to study ``equimodular'' limit periodic continued fractions, which
are limit periodic continued fractions in which the characteristic roots of the associated
$2 \times 2$ matrices are all equal in modulus. (Thus they are exactly
the class of limit periodic continued fractions of elliptic type.)
The matrix theorem in \cite{KL} can
also be used to obtain results about the boundedness of recurrence sequences.
Theorem \ref{t2} below applies to equimodular recurrences as
well.

More recent is the work of R.J. Kooman \cite{koomR,koomM,koom}.
Kooman makes a detailed study of the asymptotics of Poincar\'e type
recurrences as well as outer composition sequences of M\"obius transformations.
Following our theorem, we compare our theorem with results of Kooman.

Our focus is on the case where the characteristic roots are of
equal modulus but distinct.
Under an $l_1$ perturbation
we will show that all non-trivial solutions of such
recurrences are asymptotic to a linear recurrence with constant
coefficients. Our theorem is:

\begin{theorem}\label{t2}
Let the sequence $\{x_{n}\}_{n \geq 0}$ be defined by initial values
$x_{0}$, $\dots$, $x_{p-1}$ and by \eqref{pereq}
for $n \geq 0$. Suppose also that there are numbers $a_{0},\dots,
a_{p-1}$ such that
\begin{align*}
&\sum_{n=0}^{\infty}|a_{r}-a_{n,r}|<\infty,& &0 \leq r \leq p-1.&
\end{align*}
Put
\[
\varepsilon_n=\max_{0\le r < p}\left(\sum_{i>n}|a_r-a_{i,r}| \right).
\]
Suppose further that the roots of the characteristic equation
\begin{equation}\label{chareqa}
t^{p}-a_{p-1}t^{p-1}-a_{p-2}t^{p-2}- \dots -a_{0}=0
\end{equation}
 are distinct with values
$\alpha_{0}$, $\dots$, $\alpha_{p-1}$ of equal modulus $R>0$.
Then there exist complex
numbers $c_0, \dots , c_{p-1}$ such that
\begin{equation}\label{xrecur}
R^{-n} \left| x_{n} - \sum_{i=0}^{p-1}c_{i}\alpha_{i}^{n}\right|=O\left(\varepsilon_n\right).
\end{equation}
\end{theorem}
\begin{proof}
Assume first that all the characteristic roots have modulus $1$. Define
\begin{equation*}
M:= \left (
\begin{matrix}
&a_{p-1} &a_{p-2}& \dots & a_{1}&a_{0}\\
&1 & 0& \dots & 0& 0\\
&0 & 1& \dots & 0& 0\\
&\vdots & \vdots& \ddots & \vdots& \vdots\\
&0 & 0& \dots & 1& 0\\
\end{matrix}
\right ).
\end{equation*}
By the correspondence between polynomials and companion matrices,
the eigenvalues of $M$ are $\alpha_{0}, \dots, \alpha_{p-1}$, so that
$M$ is diagonalizable. For $n \geq 1$, define
\begin{equation*}
D_{n}:= \left (
\begin{matrix}
&a_{n-1,p-1} &a_{n-1,p-2}&
                    \dots & a_{n-1,1}&a_{n-1,0}\\
&1 & 0& \dots & 0& 0\\
&0 & 1& \dots & 0& 0\\
&\vdots & \vdots& \ddots & \vdots& \vdots\\
&0 & 0& \dots & 1& 0\\
\end{matrix}
\right ).
\end{equation*}
Thus the matrices $M$ and $D_{n}$ satisfy the conditions of
Theorem \ref{gent}. From  \eqref{pereq}
it follows that
\begin{equation*}
\left (
\begin{matrix}
&x_{n+p-1}\\
&x_{n+p-2}\\
& \vdots \\
& x_{n}
\end{matrix}
\right ) =\prod_{j=1}^{n}D_{j} \left (
\begin{matrix}
&x_{p-1}\\
&x_{p-2}\\
& \vdots \\
& x_{0}
\end{matrix}
\right ).
\end{equation*}
Let $F$ have the same meaning as in Theorem \ref{gent}. Part (i) then gives that
\begin{equation*}
\left| \left (
\begin{matrix}
&x_{n+p-1}\\
&x_{n+p-2}\\
& \vdots \\
& x_{n}
\end{matrix}
\right ) -F\, M^{n} \left (
\begin{matrix}
&x_{p-1}\\
&x_{p-2}\\
& \vdots \\
& x_{0}
\end{matrix}
\right )\right|=O\left(\varepsilon_n\right).
\end{equation*}
\eqref{xrecur} follows immediately by considering
the bottom entry on the left side. The case of modulus $R$ follows
by renormalization. This completes the
proof.
\end{proof}

The following corollary,  proved in \cite{BMcL05}, is immediate.

\begin{corollary}
Let the sequence $\{x_{n}\}_{n \geq 0}$ be defined by initial values
$x_{0}$, $\dots$, $x_{p-1}$ as well as \eqref{pereq} for $n \geq 0$.
Suppose also that there are numbers $a_{0},\dots, a_{p-1}$ such that
\begin{align*}
&\sum_{n=0}^{\infty}|a_{r}-a_{n,r}|<\infty, & & 0 \leq r \leq p-1.&
\end{align*}
Assume that  the roots of the characteristic equation
\begin{equation*}
t^{p}-a_{p-1}t^{p-1}-a_{p-2}t^{p-2}- \dots -a_{0}=0
\end{equation*}
are distinct roots of unity $\alpha_{0}$, $\dots$,
$\alpha_{p-1}$. Let $m$ be the least positive integer such that, for
all $j \in \{0,1,\dots,p-1\}$, $\alpha_{j}^{m} = 1$.
 Then, for $0 \leq j \leq m-1$, the subsequence
$\{x_{mn+j}\}_{n=0}^{\infty}$ converges. Set $l_j=\lim_{n\to\infty}
x_{nm+j}$, for integers $j\ge 0$. Then the (periodic) sequence
$\{l_j\}$ satisfies the recurrence relation
\begin{equation*}
l_{n+p}=\sum_{r=0}^{p-1}a_{r}l_{n+r},
\end{equation*}
and thus there exist constants $c_0,\cdots,c_{p-1}$ such that
\begin{equation*}
l_n=\sum_{i=0}^{p-1} c_i\alpha_i^n.
\end{equation*}
\end{corollary}

We close this section by comparing our result with those from Kooman \cite{koom}.
Proposition 1.7
from \cite{koom} appears to be most closely related to Theorem \ref{t2}.
Kooman also gives a result of Evgrafov \cite{Ev} which is also similar
to Theorem \ref{t2}:

{\it Consider the linear recurrence \eqref{pereq} where
\begin{align*}
&\sum_{n=0}^{\infty}|a_{r}-a_{n,r}|<\infty,& &0 \leq r \leq p-1.&
\end{align*}
If the characteristic polynomial \eqref{chareq} has zeros $\alpha_0,\dots,\alpha_{p-1}$ with
$0<|\alpha_0|\le\cdots\le|\alpha_{p-1}|$, then \eqref{pereq} has solutions
$u_n^{(i)}=\alpha_i^n(1+o(1))$.}

One difference between Evgrafov's result and Theorem \ref{t2}, is that
the later gives an error term. However, Evgrafov's theorem
does not require distinct characteristic roots. Kooman obtained a result
generalizing Evgrafov's, and containing an error term:

{\it (Proposition 1.7 of \cite{koom})  Let $\alpha_0,\dots,\alpha_{p-1}$
be non-zero, not necessarily distinct numbers with $|\alpha_0|\le\cdots\le|\alpha_{p-1}|$
and let $\beta:\mathbb{N}\to\mathbb{R}_{>0}$ be a function such that
$\lim_{n\to\infty}\beta(n)=0$, $\sum_{n=0}^{\infty}\beta(n)<\infty$, and
$0<\max |\alpha_i/\alpha_{i+1}|<\liminf(\beta(n+1)/\beta(n))\le 1$ where
the maximum is taken over those $i$ such that $|\alpha_i|\neq|\alpha_{i+1}|$.
Let $K_n$ be matrices with $\|K_n\|=O(\beta(n))$. The matrix recurrence
\[
(\textnormal{diag}(\alpha_0,\dots,\alpha_{p-1})+K_n)x_n=x_{n+1}\qquad (n\in\mathbb{N})
\]
has solutions $\{x_n^{(i)}\}$ with
\[
x_n^{(i)}=\alpha_i^ne_i\left(1+O\left(\sum_{h=n}^\infty \beta(h)\right)\right)
\]
for $i=0,\dots,p-1$.
}

Here $e_i$ is the $i$th canonical basis element for $\mathbb{C}^p$.
The asymptotic in Theorem \ref{t2} can be obtained from Kooman's Proposition 1.7
as a special case, but with a different error term. Kooman's error term can
be  weaker because of the assumption that
$$\liminf_{n\to\infty}\frac{\beta(n+1)}{\beta(n)}>0.$$
(Note that to obtain the result corresponding to the asymptotic in
Theorem \ref{t2} from Kooman's Proposition 1.7, one must diagonalize
the limiting companion matrix, $M$: set
$CMC^{-1}=\textnormal{diag}(\alpha_0,\dots,\alpha_{p-1})$. Then put
$K_n=CD_nC^{-1}-\textnormal{diag}(\alpha_0,\dots,\alpha_{p-1})$.
One easily checks that the sequence $K_n$ satisfies the conditions
of Kooman's proposition.)

\section{Conclusion}

We have studied convergent subsequences of approximants of complex
continued fractions and generalizations. There is an interesting
pattern of relationships between the limits and asymptotics of
subsequences and the modified approximants of the original sequence.
This suggests the general question of {\it in which other situations
do ``similar'' patterns of relationships exist? } In section 2, it was
shown that (at least some of) this behavior extends to the setting
of products of invertible elements in Banach algebras. More generally,
are there other classes of sequences that diverge by oscillation,
but for which ``nice'' asymptotics for the sequences exist, thus
enabling the computation of the sequential closure? When can
the probability density functions be computed? Even more
generally, when ``nice'' asymptotics do not exist, is the sequential
closure interesting or useful?

\medskip
\noindent
{\bf Acknowledgements.} The authors would like to thank
the following people for helpful discussions or comments: Bruce Berndt,
Daniel Grubb, Doug Hensley, Paul Levrie, Bruce Reznick, Ian Short, and Peter Waterman.

\allowdisplaybreaks{

}
\end{document}